\documentclass[a4paper,10pt,leqno]{article}
\usepackage{amsmath, amsfonts, amssymb, amsthm}
\numberwithin{equation}{section}

\begin{document}

\title{ On the modified multi-component Camassa-Holm system in higher dimensions}

\author{Kai Yan\footnote{email: kaiyan@hust.edu.cn }\\
School of Mathematics and Statistics,\\
Huazhong University of Science and Technology,\\
Wuhan 430074,  China}
\date{}
\maketitle

\begin{abstract}
This paper is devoted to the Cauchy problem for the modified multi-component Camassa-Holm system in higher dimensions.
On the one hand, we establish an almost complete local well-posedness results for the system in the framework of Besov spaces.
On the other hand, several blow-up criteria of strong solutions to the system are derived 
by using the Littlewood-Paley decomposition and the energy method.\\

\noindent {\bf Mathematics Subject Classification (2010)} 35Q35, 35Q51, 35L30\\

\noindent
\textbf{Keywords}: Multi-component Camassa-Holm system, higher dimensions, local well-posedness,
blow up, Besov spaces, Littlewood-Paley decomposition, energy method.
\end{abstract}

\section{Introduction}
\newtheorem {remark1}{Remark}[section]
\newtheorem{theorem1}{Theorem}[section]
\newtheorem{lemma1}{Lemma}[section]
\newtheorem{definition1}{Definition}[section]
\newtheorem{corollary1}{Corollary}[section]
\newtheorem{proposition1}{Proposition}[section]

In this paper, we consider the Cauchy problem for the following modified multi-component Camassa-Holm system in higher dimensions:
\begin{equation}\label{1mch}
\left\{\begin{array}{ll}
\partial_t m+ \underbrace{u\cdot\nabla m}_{\text{convection}}+ \underbrace{\nabla u^{T}\cdot m}_{\text{stretching}}
+ \underbrace{m(\text{div} u)}_{\text{expansion}}+\underbrace{\rho \nabla\bar{\rho}}_{\text{force}}=0,\\
\partial_t\rho+\text{div}(\rho u)=0,
\end{array}\right.
\end{equation}
or in components,
\begin{equation*}
\left\{\begin{array}{ll}
\frac{\partial m_i}{\partial t}+ \sum\limits_{j=1}^d u_j \frac{\partial m_i}{\partial x_j}
+ \sum\limits_{j=1}^d m_j \frac{\partial u_j}{\partial x_i}
+ m_i \sum\limits_{j=1}^d \frac{\partial u_j}{\partial x_j}+\rho \frac{\partial \bar{\rho}}{\partial x_i}=0,
\,\ \,\ i=1,2,\cdot\cdot\cdot,d,\\
\partial_t\rho+\sum\limits_{j=1}^d \frac{\partial (\rho u_j)}{\partial x_j}=0.
\end{array}\right.
\end{equation*}
Here the vector fields $u=u(t,x)$ and $m=m(t,x)$ are defined from $\mathbb{R^+}\times \mathbb{R}^d$
(or $\mathbb{R^+}\times \mathbb{T}^d$) to $\mathbb{R}^d$ such that $m=(I-\Delta)u$,
the scalar functions $\rho=\rho(t,x)$ and $\bar{\rho}=\bar{\rho}(t,x)$ are defined from $\mathbb{R^+}\times \mathbb{R}^d$
(or $\mathbb{R^+}\times \mathbb{T}^d$) to $\mathbb{R}$ such that $\rho=(I-\Delta)(\bar{\rho}-\bar{\rho}_0)$,
and the torus $\mathbb{T}^d \triangleq \mathbb{R}^d /\mathbb{Z}^d$.

As a set of semidirect-product Euler-Poincar\'{e} equations, the system (\ref{1mch}) was proposed in \cite{Holm}
and shown that the last four terms in the first equation of system (\ref{1mch}) model convection, stretching, expansion and force of a fluid
with velocity $u$, momentum $m$, density $\rho$ and averaged density $\bar{\rho}$, respectively.
Moreover, the system (\ref{1mch}) possess $\delta$ function-like singular solutions in both $m$ and $\rho$, which emerge from smooth initial conditions \cite{Holm}.

For $\rho\equiv 0$ and $d=1$, system (\ref{1mch}) becomes the celebrated Camassa-Holm equation (CH):
\begin{eqnarray*}
m_t+ u m_x+ 2 u_x m=0,\,\ \,\ m=u-u_{xx},
\end{eqnarray*}
which models the unidirectional propagation of shallow water waves over a flat bottom \cite{C-H}.
CH is also a model for the propagation of axially symmetric waves in hyper-elastic rods \cite{Dai}.
It has a bi-Hamiltonian structure and is completely integrable \cite{C-H}.
Its solitary waves are peaked solitons (peakons) \cite{C-H-H,Con-E}, and they are orbitally stable \cite{C-S,C-S1}.
It is noted that the peakons replicate a feature that is characteristic for the waves of great height
-- waves of the largest amplitude that are exact traveling wave solutions of the governing equations for irrotational water waves, cf. \cite{Cinvent,C-Eann}.
The Cauchy problem and initial boundary value problem for CH have been studied extensively \cite{B-C1,B-C2,C-Ep,C-Ec,Dan1,Dan2,E-Y1,HM1,R}.
It has been shown that this equation is locally well-posed \cite{C-Ep,C-Ec,Dan1,Dan2,R}.
Moreover, it has both global strong solutions \cite{Cf,C-Ep,C-Ec} and blow-up solutions within finite time \cite{Cf,C-Ep,C-E,C-Ec}.
In addition, it possess global weak solutions, see the discussions in \cite{B-C1,B-C2,X-Z}.
It is worthy to point out the advantage of CH in comparison with the KdV equation lies in the fact that
CH has peakons and models wave breaking \cite{C-H-H,C-E} (namely, the wave remains bounded while its slope becomes unbounded in finite time \cite{Wh}).

For $d=1$, system (\ref{1mch}) becomes the following modified two-component Camassa-Holm system (M2CH):
\begin{equation*}
\left\{\begin{array}{ll}
m_{t}+u m_x+2u_x m +\rho(1-\partial_x^2)^{-1}\rho_x=0,\,\ \,\ m=u-u_{xx},\\
\rho_{t}+(u\rho)_x=0,
\end{array}\right.
\end{equation*}
which was firstly proposed in \cite{Holm} and proved that it allows singular solutions in both variables $m$ and $\rho$, not just the fluid momentum.
The Cauchy problem and initial boundary value problem for (M2CH) have been investigated in many works,
see the discussions in \cite{Guan1,Guan2,Tan,Yan1,Yan3}.

For $\rho\equiv 0$, system (\ref{1mch}) reduces the higher dimensional Camassa-Holm equations as follows:
\begin{eqnarray}\label{1-d-dim-CH}
m_t+ u\cdot\nabla m+ \nabla u^{T}\cdot m+ m(\text{div} u)=0,\,\ \,\ m=(I-\Delta)u,
\end{eqnarray}
which was proposed exactly in the way that
a class of its singular solutions generalize the peakon solutions of the CH equation to higher spatial dimensions \cite{Holm03}.
It was also studied as Euler-Poincar\'{e} equations associated with the diffeomorphism group in \cite{Holm05}.
The local well-posedness in Sobolev spaces, blow up criteria, global and blow-up solutions of the Cauchy problem for Eqs. (\ref{1-d-dim-CH})
has been discussed in \cite{Chae-CMP,Duan-JFA,Li-ARMA,Yan-Yin-DCDS-2015}.

Now, let $\gamma\triangleq \bar{\rho}-\bar{\rho}_0$.
Then the Cauchy problem for system (\ref{1mch}) can be rewritten to the nonlocal form as follows (see Appendix for the details):
\begin{equation}\label{1mchnew}
\left\{\begin{array}{ll}
\partial_t u+u\cdot \nabla u= F_1(u,\gamma), &(t, x)\in \mathbb{R}\times \mathbb{R}^d,\\
\partial_t\gamma+u\cdot \nabla \gamma= F_2(u,\gamma), &(t, x)\in \mathbb{R}\times \mathbb{R}^d,\\
u(0,x)= u_{0}(x), &x\in \mathbb{R}^d,\\
\gamma(0,x)= \gamma_{0}(x), &x\in \mathbb{R}^d,
\end{array}\right.
\end{equation}
where
\begin{eqnarray}\label{1-F1}
\quad\quad\quad F_1(u,\gamma)&\triangleq &
-(I-\Delta)^{-1}\text{div}\left(\nabla u(\nabla u+\nabla u^T)-\nabla u^T \nabla u-\nabla u (\text{div} u)\right)\\ \nonumber
&&-(I-\Delta)^{-1}\text{div}\left(\frac{1}{2}(|\nabla u|^2+\gamma^2+|\nabla \gamma|^2)I-\nabla \gamma^T \nabla \gamma\right)\\ \nonumber
&&-(I-\Delta)^{-1}\left(u(\text{div} u)+u\cdot \nabla u^T\right),
\end{eqnarray}
and
\begin{eqnarray}\label{1-F2}
F_2(u,\gamma)&\triangleq&
-(I-\Delta)^{-1}\text{div}\left( \nabla \gamma \nabla u+(\nabla \gamma)\cdot \nabla u-\nabla \gamma(\text{div} u) \right)\\ \nonumber
&&-(I-\Delta)^{-1}\left(\gamma(\text{div} u)\right).
\end{eqnarray}

To our best knowledge, the Cauchy problem for system (\ref{1mch}) or the system (\ref{1mchnew}) has not been discussed yet.
It is noted that, unlike the above two-component system (M2CH) in one dimension and the single equation (CH) or equations (\ref{1-d-dim-CH}),
the present considered system is a multi-component transport equations in higher dimensions and no more regularity is available from it.
Moreover, the system (\ref{1mchnew}) is coupled with the vector field $u$ and the scalar function $\gamma$ so that
we have to deal with the mutual effect between them, for which more delicate nonlinear estimates are required in this paper.
The purpose of this paper is to establish the local well-posedness for system (\ref{1mchnew})
and derive some blow-up criteria of strong solutions to the system in the framework of Besov spaces.
Since our obtained results can be easily carried out to the periodic case
and to the homogeneous Besov spaces, we shall always assume that the space variables belong to the whole $\mathbb{R}^d$
and restrict our attention to nonhomogeneous Besov spaces.

For this, we introduce some notations. Let $s\in\mathbb{R}$, $1\leq p,r\leq\infty$.
The nonhomogeneous Besov space $B^s_{p,r}(\mathbb{R}^d)$ ($B^s_{p,r}$ for short) is defined by
$$B^s_{p,r}(\mathbb{R}^d)\triangleq \{f \in \mathcal{S'}(\mathbb{R}^d):
||f||_{B^s_{p,r}(\mathbb{R}^d)}\triangleq ||(2^{qs}||\Delta_q{f}||_{L^p(\mathbb{R}^d)})_{q\geq-1}||_{l^r}< \infty\},$$
where $\Delta_q$ is the Littlewood-Paley decomposition operator \cite{BCD} .
If $s=\infty$, then $B^\infty_{p,r}(\mathbb{R}^d)\triangleq \bigcap\limits_{s\in\mathbb{R}}B^s_{p,r}(\mathbb{R}^d)$.
Moreover, define
$$E^s_{p,r}(T)\triangleq C([0,T]; B^s_{p,r}(\mathbb{R}^d))\cap C^1{([0,T]; B^{s-1}_{p,r}(\mathbb{R}^d))},\quad \text{if}\,\ r<\infty,$$
and
$$E^s_{p,\infty}(T)\triangleq L^{\infty}(0,T; B^s_{p,\infty}(\mathbb{R}^d))\cap Lip\,(0,T; B^{s-1}_{p,\infty}(\mathbb{R}^d))$$
for some $T>0$.
In addition, if $u=(u_1,u_2,\cdot\cdot\cdot,u_d)$ is a vector field,
then for simplicity, we always write
$u\in B^s_{p,r}(\mathbb{R}^d)$ and $\nabla u\in B^s_{p,r}(\mathbb{R}^d)$
standing for $u\in (B^s_{p,r}(\mathbb{R}^d))^d$ and $\nabla u\in (B^s_{p,r}(\mathbb{R}^d))^{d^2}$, respectively,
if there is no ambiguity. And the corresponding norms notation should be understood in the same way.

In the present paper, we first obtain the following local well-posedness results in the supercritical and critical Besov spaces (Theorem 1.1 and Theorem 1.2), respectively:
\begin{theorem1}
Let $d\in \mathbb{N_+}$, $1\leq p,r\leq \infty$ and $s>\max (1+\frac{d}{p},\frac{3}{2})$.
Suppose that $(u_0, \gamma_0)\in B^s_{p,r}(\mathbb{R}^d)\times B^s_{p,r}(\mathbb{R}^d)$.
Then there exists a time $T>0$
such that $(u,\gamma)\in E^s_{p,r}(T)\times E^s_{p,r}(T)$ is the unique solution to system (\ref{1mchnew}),
and the solution depends continuously on the initial data,
that is, the mapping $(u_0,\gamma_0)\mapsto (u,\gamma)$ is continuous from
$B^s_{p,r}(\mathbb{R}^d)\times B^s_{p,r}(\mathbb{R}^d)$ into
$$ C([0,T]; B^{s'}_{p,r}(\mathbb{R}^d)\times B^{s'}_{p,r}(\mathbb{R}^d))
\cap C^1{([0,T]; B^{s'-1}_{p,r}(\mathbb{R}^d)\times B^{s'-1}_{p,r}(\mathbb{R}^d))}$$
for all $s'<s$ if $r=\infty$, and $s'=s$ otherwise.
\end{theorem1}

Note that for any $s\in\mathbb{R}$, $H^s(\mathbb{R}^d)= B^s_{2,2}(\mathbb{R}^d)$.
Then we instantaneously get the following local well-posedness result in Sobolev spaces.
\begin{corollary1}
Let $(u_0, \gamma_0)$ be in $H^s(\mathbb{R}^d)\times H^s(\mathbb{R}^d)$ with $s>1+\frac{d}{2}$.
Then there exist a time $T>0$ and a unique solution $(u,\gamma)$ to system (\ref{1mchnew}) such that
$(u,\gamma)\in C([0,T]; H^s(\mathbb{R}^d)\times H^s(\mathbb{R}^d))
\cap C^1{([0,T]; H^{s-1}(\mathbb{R}^d)\times H^{s-1}(\mathbb{R}^d))}$.
Moreover, the mapping $(u_0,\gamma_0)\mapsto (u,\gamma):$
$$ H^s\times  H^s \!\rightarrow\!  C([0,T]; H^s\times  H^s)\cap C^1{([0,T]; H^{s-1}\times H^{s-1})}$$
is continuous.
\end{corollary1}

\begin{theorem1}
Let $d\in \mathbb{N_+}$ and  $1\leq p\leq 2d$.
Suppose that $(u_0,\gamma_0)\in B^{1+\frac{d}{p}}_{p,1}(\mathbb{R}^d)\times B^{1+\frac{d}{p}}_{p,1}(\mathbb{R}^d)$.
Then there exist a time $T=T(||(u_0,\gamma_0)||_{B^{1+\frac{d}{p}}_{p,1}\times B^{1+\frac{d}{p}}_{p,1}})>0$
and a unique solution $(u,\gamma)$ to system (\ref{1mchnew}) such that
$$(u,\gamma)\in C([0,T]; B^{1+\frac{d}{p}}_{p,1}(\mathbb{R}^d)\times B^{1+\frac{d}{p}}_{p,1}(\mathbb{R}^d))
\cap C^1{([0,T]; B^{\frac{d}{p}}_{p,1}(\mathbb{R}^d)\times B^{\frac{d}{p}}_{p,1}(\mathbb{R}^d))}.$$
Moreover, the solution depends continuously on the initial data,
that is, the mapping $(u_0,\gamma_0)\mapsto (u,\gamma)$ is continuous from
$B^{1+\frac{d}{p}}_{p,1}(\mathbb{R}^d)\times B^{1+\frac{d}{p}}_{p,1}(\mathbb{R}^d)$ into
$$C([0,T]; B^{1+\frac{d}{p}}_{p,1}(\mathbb{R}^d)\times B^{1+\frac{d}{p}}_{p,1}(\mathbb{R}^d))
\cap C^1{([0,T]; B^{\frac{d}{p}}_{p,1}(\mathbb{R}^d)\times B^{\frac{d}{p}}_{p,1}(\mathbb{R}^d))}.$$
\end{theorem1}

\begin{remark1}
(1) Theorem 1.1 and Theorem 1.2 cover and extend the corresponding results in \cite{Chae-CMP,Dan2,Guan1,Yan-Yin-DCDS-2015}.
Moreover, Corollary 1.1 in the case of $d=2$ improves the related result in \cite{Kman},
where the periodic 2D Camassa-Holm equations is proved locally well-posed as the initial data $u_0\in H^s(\mathbb{T}^2)$ with $s>3$ by using a geometric approach.\\
(2) Note that for any $s>1+\frac{d}{p}$, $B^{s}_{p,r}(\mathbb{R}^d)\hookrightarrow B^{1+\frac{d}{p}}_{p,1}(\mathbb{R}^d)$.
Theorem 1.2 improves the corresponding result in Theorem 1.1 when $1\leq p\leq 2d$.
However, except the existence of the solutions, the question of uniqueness and continuity with respect to the initial data
$(u_0,\gamma_0)\in B^{1+\frac{d}{p}}_{p,1}(\mathbb{R}^d)\times B^{1+\frac{d}{p}}_{p,1}(\mathbb{R}^d)$ as $2d<p<\infty$ still remains unknown.\\
(3) It is well known that for any $s'<1+\frac{d}{2}<s$, the following embedding relations
$$H^s(\mathbb{R}^d)\hookrightarrow B^{1+\frac{d}{2}}_{2,1}(\mathbb{R}^d)\hookrightarrow
H^{1+\frac{d}{2}}(\mathbb{R}^d)\hookrightarrow B^{1+\frac{d}{2}}_{2,\infty}(\mathbb{R}^d)\hookrightarrow H^{s'}(\mathbb{R}^d)$$
hold true, which shows that $H^s(\mathbb{R}^d)$ and $B^s_{2,r}(\mathbb{R}^d)$ are quite close.
Corollary 1.1 and Theorem 1.2 ensure the local well-posedness for system (\ref{1mchnew})
as initial data $(u_0,\gamma_0)$ belongs to $H^s(\mathbb{R}^d)\times H^s(\mathbb{R}^d)$ or $B^{1+\frac{d}{2}}_{2,1}(\mathbb{R}^d)\times B^{1+\frac{d}{2}}_{2,1}(\mathbb{R}^d)$.
However, whether the system (\ref{1mchnew}) is locally well-posed or not when $(u_0,\gamma_0)\in H^{1+\frac{d}{2}}(\mathbb{R}^d)\times H^{1+\frac{d}{2}}(\mathbb{R}^d)$ is an open problem so far.
While it is noted to point out that, as a special case of system (\ref{1mchnew}), the 1D Camassa-Holm equation with initial data
$u_0\in B^{\frac{3}{2}}_{2,\infty}(\mathbb{R})$ or $u_0\in H^{s'}(\mathbb{R})$ with $s'<\frac{3}{2}$
is not locally well-posed in the sense that the solutions do not depend uniformly continuously on the initial data, cf. \cite{Dan2,HM1}.
So, in this context, the system (\ref{1mchnew}) is ill-posed in the subcritical Besov spaces (i.e. the regularity index $s<1+\frac{d}{p}$),
and the results in Theorem 1.1 and Theorem 1.2 are sharp.
Overall, we here give an almost complete local well-posedness results for system (\ref{1mchnew}) in the framework of Besov spaces.
\end{remark1}

Next, we prove three blow-up criteria (Theorems 1.3-1.5) of the strong solutions to system (\ref{1mchnew}) as follows.
\begin{theorem1}
Suppose that $d\in\mathbb{N_+}$ and $1\leq p,r\leq \infty$.
Let $(u_0, \gamma_0)\in B^s_{p,r}(\mathbb{R}^d)\times B^s_{p,r}(\mathbb{R}^d)$ with
$s>\max(1+\frac{d}{p},\frac 3 2)$ (or $s=1+\frac{d}{p}$ with $r=1$ and $1\leq p\leq 2d$),
and $(u,\gamma)$ be the corresponding solution to system (\ref{1mchnew}).
Then the solution blows up in finite time (i.e. the lifespan of solution $T^{\star}<\infty$ ) if and only if
$$\int_0^{T^{\star}}
\left(||u(\tau)||_{L^\infty(\mathbb{R}^d)}+||\nabla u(\tau)||_{L^\infty(\mathbb{R}^d)}+||\gamma(\tau)||_{L^\infty(\mathbb{R}^d)}+||\nabla\gamma(\tau)||_{L^\infty(\mathbb{R}^d)}\right) d\tau= \infty.$$
\end{theorem1}

\begin{remark1}
It is not difficult to check that
$$H(d)\triangleq \int_{\mathbb{R}^d} (|u|^2+|\nabla u|^2+\gamma^2+|\nabla \gamma|^2) d x$$
is a conservation law of system (\ref{1mchnew}).
The Sobolev embedding theorem implies that
$H(d)$ cannot be utilized to control $||u||_{L^\infty(\mathbb{R}^d)}+||\gamma||_{L^\infty(\mathbb{R}^d)}$ unless $d=1$.
So, Theorem 1.3 in the case of $d=1$ can be improved as
$$T^{\star}<\infty \Leftrightarrow \int_0^{T^{\star}} \left(||\nabla u(\tau)||_{L^\infty(\mathbb{R})}+||\nabla\gamma(\tau)||_{L^\infty(\mathbb{R})}\right) d\tau= \infty.$$
While beyond one dimension, the above approach is no longer valid.
Fortunately, by further exploring the structures of system (\ref{1mchnew}),
one can obtain a more precise blow-up criterion for arbitrary dimensions (see Theorem 1.4 below), which depends only on $\nabla u$ .
\end{remark1}

\begin{theorem1}
Under the assumptions in Theorem 1.3, the solution $(u,\gamma)$ blows up in finite time (i.e. the lifespan of solution $T^{\star}<\infty$ ) if and only if
$$\int_0^{T^{\star}} ||\nabla u(\tau)||_{L^\infty(\mathbb{R}^d)}d\tau= \infty.$$
\end{theorem1}

\begin{remark1} \label{1-maximal-time}
The maximal existence time $T$ in Theorem 1.1 can be chosen independent of the regularity index $s$.
Indeed, let $(u_0,\gamma_0)\in B^s_{p,r}\times B^s_{p,r}$ with $s>\max(1+\frac{d}{p},\frac 3 2)$ and some $s'\in (\max(1+\frac{d}{p},\frac 3 2),s)$.
Then Theorem 1.1 ensures that there exists a unique $B^s_{p,r}\times B^s_{p,r}$ (resp., $B^{s'}_{p,r}\times B^{s'}_{p,r}$)
solution $(u_s,\gamma_s)$ (resp., $(u_{s'},\gamma_{s'})$) to system (\ref{1mchnew}) with the maximal existence time $T_s$ (resp., $T_{s'}$).
Since $B^s_{p,r}\hookrightarrow B^{s'}_{p,r}$, it follows from the uniqueness that $T_s\leq T_{s'}$ and $u_s\equiv u_{s'}$ on $[0, T_s)$.
On the other hand, if we suppose that $T_s < T_{s'}\leq \infty$,
then $u_{s'}\in C([0,T_s]; B^{s'}_{p,r})$. Hence, $\nabla u_s\in L^1(0,T_s; L^\infty)$, which is a contradiction to Theorem 1.4.
Therefore, $T_s = T_{s'}$.
Likewise, denote $T_{\text{critical}}$ and $T_{r=1}$ by the maximal existence times in Theorem 1.2 and in Theorem 1.1 with $r=1$, respectively.
Then $T_{\text{critical}}=T_{r=1}$.
\end{remark1}

Notice that for any $s>\frac{d}{p}$ (or $s= \frac{d}{p}$ and $r=1$), we have
$$B^s_{p,r}(\mathbb{R}^d)\hookrightarrow L^\infty(\mathbb{R}^d)\hookrightarrow BMO(\mathbb{R}^d)\hookrightarrow B^0_{\infty,\infty}(\mathbb{R}^d).$$
Finally, we derive a blow-up criterion in terms of the $B^0_{\infty,\infty}(\mathbb{R}^d)$ norm.
While compared to the result in Theorem 1.4, the cost we pay here is that both $\nabla u$ and $\nabla \gamma$ will be involved.
\begin{theorem1}
Under the assumptions in Theorem 1.3, the solution $(u,\gamma)$ blows up in finite time (i.e. the lifespan of solution $T^{\star}<\infty$ ) if and only if
$$\int_0^{T^{\star}} \left(||\nabla u(\tau)||_{B^0_{\infty,\infty}(\mathbb{R}^d)}+||\nabla\gamma(\tau)||_{B^0_{\infty,\infty}(\mathbb{R}^d)}\right)d\tau= \infty.$$
\end{theorem1}

The rest of our paper is organized as follows.
In Section 2, we recall some fine properties of Besov spaces and the transport equations theory.
In Section 3, we prove Theorem 1.1 to establish the local well-posedness of system (\ref{1mchnew}) in supercritical Besov spaces.
In Section 4, we prove Theorem 1.2 to establish the local well-posedness of system (\ref{1mchnew}) in critical Besov spaces.
In Section 5, we derive the blow-up criteria of strong solutions to system (\ref{1mchnew}) by showing Theorems 1.3-1.5.
Section 6 is devoted to an Appendix.

\section{Preliminaries}
\newtheorem {remark2}{Remark}[section]
\newtheorem{theorem2}{Theorem}[section]
\newtheorem{lemma2}{Lemma}[section]
\newtheorem{definition2}{Definition}[section]
\newtheorem{proposition2}{Proposition}[section]
\newtheorem{corollary2}{Corollary}[section]

In this section, we recall some fine properties of Besov spaces and the transport equations theory, which are frequently used in the whole paper.
\begin{lemma2}
\cite{BCD}
(i) Complex interpolation: $\forall\, f\in B^{s_1}_{p,r}(\mathbb{R}^d)\cap B^{s_2}_{p,r}(\mathbb{R}^d)$,
\begin{eqnarray} \label{2interpolation}
||f||_{B^{\theta {s_1}+(1-\theta){s_2}}_{p,r}(\mathbb{R}^d)}\leq ||f||^{\theta}_{B^{s_1}_{p,r}(\mathbb{R}^d)}||f||^{1-\theta}_{B^{s_2}_{p,r}(\mathbb{R}^d)},
\,\ \,\ \theta \in[0,1].
\end{eqnarray}
(ii) Logarithmic type interpolation inequality:
there exists a positive constant $C$ such that for all
$s\in\mathbb{R}$, $\varepsilon>0$ and $1\leq p\leq\infty$, we have
\begin{eqnarray} \label{2log-interpolation}
||f||_{B^{s}_{p,1}(\mathbb{R}^d)}
\leq C\frac{1+\varepsilon}{\varepsilon}||f||_{B^{s}_{p,\infty}(\mathbb{R}^d)}
\ln \left(e+\frac{||f||_{B^{s+\varepsilon}_{p,\infty}(\mathbb{R}^d)}}
{||f||_{B^{s}_{p,\infty}(\mathbb{R}^d)}}\right).
\end{eqnarray}
\end{lemma2}

\begin{corollary2}
There exists a positive constant $c$ such that for any $q>d$,
\begin{eqnarray} \label{2new-log-interpolation}
\quad\quad ||f||_{L^\infty(\mathbb{R}^d)}
\leq c \frac{2q-d}{q-d}
\left(1+||f||_{B^{0}_{\infty,\infty}(\mathbb{R}^d)} \ln(e+||f||_{W^{1,q}(\mathbb{R}^d)})\right).
\end{eqnarray}
\end{corollary2}

\begin{proof}
Note that $B^0_{\infty,1}(\mathbb{R}^d)\hookrightarrow L^\infty(\mathbb{R}^d)$ and
$W^{1,q}(\mathbb{R}^d)\hookrightarrow B^{1-\frac{d}{q}}_{\infty,\infty}(\mathbb{R}^d)$.
Set $s=0$, $p=\infty$ and $\varepsilon=1-\frac{d}{q}>0$ in (\ref{2log-interpolation}).
Then one infers that
\begin{eqnarray*}
||f||_{L^\infty(\mathbb{R}^d)}
&\leq& c ||f||_{B^0_{\infty,1}(\mathbb{R}^d)}\\
&\leq& c \frac{2q-d}{q-d} ||f||_{B^0_{\infty,\infty}(\mathbb{R}^d)}
\ln\left(e+\frac{||f||_{B^{1-\frac{d}{q}}_{\infty,\infty}(\mathbb{R}^d)}} {||f||_{B^0_{\infty,\infty}(\mathbb{R}^d)}}\right)\\
&\leq& c \frac{2q-d}{q-d} ||f||_{B^0_{\infty,\infty}(\mathbb{R}^d)}
\ln\left(e+\frac{||f||_{W^{1,q}(\mathbb{R}^d)}} {||f||_{B^0_{\infty,\infty}(\mathbb{R}^d)}}\right)\\
&\leq& c \frac{2q-d}{q-d}
\left(1+||f||_{B^{0}_{\infty,\infty}(\mathbb{R}^d)} \ln(e+||f||_{W^{1,q}(\mathbb{R}^d)})\right),
\end{eqnarray*}
where $c$ is independent of $q$.
\end{proof}

\begin{lemma2} \label{2-multiplier}
\cite{BCD}
Let $m\in\mathbb{R}$ and $f$ be an $S^m$-multiplier.
That is, $f: \mathbb{R}^d \rightarrow \mathbb{R}$ is smooth and satisfies that
for any $\alpha\in \mathbb{N}^d$, there is a constant $C_{\alpha}>0$ such that
$$|\partial^{\alpha}{f(\xi)}|\leq C_{\alpha}(1+|\xi|)^{m-|\alpha|},\,\ \,\ \forall\,\ \xi\in\mathbb{R}^d.$$
Set $f(D)\triangleq\mathcal{F}^{-1} f \mathcal{F}\in Op(S^m)$.
Then the operator $f(D)$ is continuous from $B^s_{p,r}(\mathbb{R}^d)$ to $B^{s-m}_{p,r}(\mathbb{R}^d)$.
\end{lemma2}

\begin{lemma2}\cite{BCD,Yan-Yin-DCDS-2015}
($d$-dimensional Morse-type estimates) Let $d\in\mathbb{N_+}$. Then\\
(i) For any $s>0$ and $1\leq p,r\leq\infty$, there exists $C=C(d,s,p,r)>0$ such that
\begin{eqnarray}\label{2morse1}
\quad\quad\quad ||f g||_{B^s_{p,r}(\mathbb{R}^d)}
\leq C\left(||f||_{L^\infty(\mathbb{R}^d)} ||g||_{B^s_{p,r}(\mathbb{R}^d)}
+||g||_{L^\infty(\mathbb{R}^d)} ||f||_{B^s_{p,r}(\mathbb{R}^d)}\right).
\end{eqnarray}
(ii) If $1\leq p,r\leq\infty$, $s_1\leq \frac{d}{p}<s_2$ ($s_2\geq {d\over p}$ if $r=1$) and $s_1+s_2>0$,
then there exists $C=C(s_1,s_2,d,p,r)>0$ such that
\begin{eqnarray}\label{2morse2}
||f g||_{B^{s_1}_{p,r}(\mathbb{R}^d)}
\leq C||f||_{B^{s_1}_{p,r}(\mathbb{R}^d)} ||g||_{B^{s_2}_{p,r}(\mathbb{R}^d)}.
\end{eqnarray}
(iii) A critical Morse-type estimate \cite{Yan-Yin-DCDS-2015}:\\
If $1\leq p\leq 2d$, then there exists $C=C(d,p)>0$ such that
\begin{eqnarray}\label{2morsecritical}
||f g||_{B^{\frac{d}{p}-1}_{p,\infty}(\mathbb{R}^d)}
\leq C ||f||_{B^{\frac{d}{p}-1}_{p,1}(\mathbb{R}^d)}
||g||_{B^{\frac{d}{p}}_{p,\infty}(\mathbb{R}^d)\cap L^\infty(\mathbb{R}^d)}.
\end{eqnarray}
\end{lemma2}

\begin{lemma2}\label{2commutator}
\cite{BCD} (Commutator estimates)
Assume $d\in\mathbb{N_+}$, $1\leq p,r\leq \infty$ and $s>0.$ Let $v$ be a vector field over $\mathbb{R}^d$.
Then there exists a constant $C=C(d,s,p)>0$ such that
\begin{eqnarray*}
&&||(2^{qs}||[v, \Delta_q]\cdot\nabla f||_{L^p(\mathbb{R}^d)})_{q\geq-1}||_{l^r}\\ \nonumber
&\leq& C\left(||\nabla v||_{L^\infty(\mathbb{R}^d)} ||f||_{B^s_{p,r}(\mathbb{R}^d)}
+||\nabla f||_{L^\infty(\mathbb{R}^d)} ||\nabla v||_{B^{s-1}_{p,r}(\mathbb{R}^d)}\right),
\end{eqnarray*}
where $[A,B]\triangleq AB-BA$ is the commutator for two operators $A$ and $B$.
\end{lemma2}

Next, we state a priori estimates for the transport equations in Besov spaces as follows.
\begin{lemma2}\label{2transport-estimates}
\cite{BCD}
Let $d\in\mathbb{N_+}$, $1\leq p,r\leq \infty$ and $s>-\min ({d\over p}, 1-{d\over p}).$
Assume that $f_0\in B^s_{p,r}(\mathbb{R}^d)$, $F\in L^1(0,T; B^s_{p,r}(\mathbb{R}^d))$,
and $\nabla v$ belongs to $L^1(0,T; B^{s-1}_{p,r}(\mathbb{R}^d))$ if $s> 1+{d\over p}$,
or to $L^1(0,T; B^{d\over p}_{p,r}\cap L^\infty(\mathbb{R}^d))$ otherwise.
If $f\in L^\infty(0,T; B^s_{p,r}(\mathbb{R}^d)) \cap C([0,T]; \mathcal{S'}(\mathbb{R}^d))$
solves the following transport equations:
\[(TE)\left\{\begin{array}{l}
\partial_t f+v\cdot\nabla f=F,\\
f|_{t=0} =f_0,
\end{array}\right.\]
then there exists a constant $C=C(d,s,p,r)>0$ such that\\
(i) If $r=1$ or $s\neq 1+{d\over p}$,
\begin{equation*}
||f(t)||_{B^s_{p,r}(\mathbb{R}^d)}\leq
||f_0||_{B^s_{p,r}(\mathbb{R}^d)}\,+\, \int_0^t ||F(\tau)||_{B^s_{p,r}(\mathbb{R}^d)}d\tau
\,+\, C\int_0^t V'(\tau)||f(\tau)||_{B^s_{p,r}(\mathbb{R}^d)} d\tau
\end{equation*}
or hence,
\begin{equation*}
||f(t)||_{B^s_{p,r}(\mathbb{R}^d)}\leq
e^{CV(t)} \left(||f_0||_{B^s_{p,r}(\mathbb{R}^d)}\,+\,
\int_0^t e^{-CV(\tau)} ||F(\tau)||_{B^s_{p,r}(\mathbb{R}^d)}d\tau\right)
\end{equation*}
with $V(t)=\int_0^t ||\nabla v(\tau)||_{B^{d\over p}_{p,r}\cap L^\infty}d\tau$
if $s<1+{d\over p}$, and $V(t)=\int_0^t ||\nabla v(\tau)||_{B^{s-1}_{p,r}}d\tau$ else.\\
(ii) If $r<\infty$, then $f\in C([0,T]; B^s_{p,r}(\mathbb{R}^d))$.
If $r=\infty$, then $f\in C([0,T]; B^{s'}_{p,1})$ for all $s'<s$.
\end{lemma2}


Finally, we need the following Osgood lemma which is a generalization of the Gronwall inequality.
\begin{lemma2}\label{2osgood}
\cite{BCD} Let $f$ be a positive measurable function, $\lambda$ a positive locally integrable function and $\mu$ a positive increasing continuous function. If
\begin{eqnarray*}
f(t)\leq \alpha+\int_{t_0}^t \lambda(s)\mu(f(s))d s,\quad \text{for}\,\ \alpha >0,
\end{eqnarray*}
then
\begin{eqnarray*}
W(f(t))\leq W(\alpha)+\int_{t_0}^t \lambda(s)d s
\end{eqnarray*}
with $W(x)\triangleq \int_a^x \frac{d r}{\mu(r)}$ for some $a>0$.
\end{lemma2}

\section{Local well-posedness in supercritical Besov spaces}
\newtheorem{remark3}{Remark}[section]
\newtheorem{theorem3}{Theorem}[section]
\newtheorem{lemma3}{Lemma}[section]
\newtheorem{corollary3}{Corollary}[section]
\newtheorem{proposition3}{Proposition}[section]

In this section, we will establish the local well-posedness of system (\ref{1mchnew}) in the supercritical Besov spaces
by using the Friedrichs regularization method and transport equations theory.

In order to prove Theorem 1.1, we first establish a priori estimates of the solutions,
which implies uniqueness and continuity with respect to the initial data in some sense.
\begin{lemma3}
Let $d$, $s$, $p$ and $r$ be in the statement of Theorem 1.1. Suppose that
$(u,\gamma), (v,\eta)\in L^{\infty}(0,T; B^s_{p,r}(\mathbb{R}^d)\times B^s_{p,r}(\mathbb{R}^d))\cap C([0,T];\mathcal{S'}(\mathbb{R}^d)\times \mathcal{S'}(\mathbb{R}^d))$
are two solutions to system (\ref{1mchnew}) with the initial data $(u_0,\gamma_0), (v_0,\eta_0)\in B^s_{p,r}(\mathbb{R}^d)\times B^s_{p,r}(\mathbb{R}^d)$, respectively.
Set $(a,b)\triangleq (v-u,\eta-\gamma)$ and $(a_0,b_0)\triangleq (v_0-u_0,\eta_0-\gamma_0)$.
Then for all $t\in[0,T]$, we have\\
(i) if $s>\max (1+\frac{d}{p},\frac 3 2)$ and $s\neq 2+\frac{d}{p}$, or $r=1$,
\begin{eqnarray}\label{3-lemma1-i}
&&||a(t)||_{B^{s-1}_{p,r}}+||b(t)||_{B^{s-1}_{p,r}}\\ \nonumber
&\leq& (||a_0||_{B^{s-1}_{p,r}}+||b_0||_{B^{s-1}_{p,r}})
e^{C\int_0^t (||u(\tau)||_{B^s_{p,r}}+||\gamma(\tau)||_{B^s_{p,r}}+||v(\tau)||_{B^s_{p,r}}+||\eta(\tau)||_{B^s_{p,r}})d\tau}\\ \nonumber
&\triangleq& M(t;s-1);
\end{eqnarray}
(ii) if $s= 2+\frac{d}{p}$ and $r\neq 1$,
\begin{eqnarray} \label{3-lemma1-ii}
&& ||a(t)||_{B^{s-1}_{p,r}}+||b(t)||_{B^{s-1}_{p,r}}\\ \nonumber
&\leq& C M^{\theta}(t;s-1) (||u(t)||_{B^s_{p,r}}+||\gamma(t)||_{B^s_{p,r}}+||v(t)||_{B^s_{p,r}}+||\eta(t)||_{B^s_{p,r}})^{1-\theta},
\end{eqnarray}
where $\theta\in(0,1)$ and $C=C(d,s,p,r)>0$.
\end{lemma3}

\begin{proof}
Apparently, $(a,b)\in L^{\infty}(0,T; B^s_{p,r}\times B^s_{p,r})\cap C([0,T];\mathcal{S'}\times \mathcal{S'})$
solves the following Cauchy problem of the transport equations:
\begin{equation} \label{3diffence-eqs}
\left\{\begin{array}{ll}
\partial_t a+u\cdot\nabla a=R_1(t,x),\\
\partial_t b+u\cdot\nabla b=R_2(t,x),\\
a|_{t=0}=a_0(x),\\
b|_{t=0}=b_0(x),
\end{array}\right.
\end{equation}
where
\begin{eqnarray*}
R_1(t,x)&\triangleq &
-a\cdot\nabla v-(I-\Delta)^{-1}\left(a(\text{div} v)+u(\text{div} a)+a\cdot \nabla v^T+u\cdot \nabla a^T\right)\\
&&-(I-\Delta)^{-1}\text{div}\left(\nabla a(\nabla v+\nabla v^T)+(\nabla u-\nabla u^T)\nabla a+\nabla u \nabla a^T-\nabla a^T \nabla v\right)\\
&&-(I-\Delta)^{-1}\text{div}\left(-\nabla a (\text{div} v)-\nabla u (\text{div} a)-\nabla \gamma^T \nabla b-\nabla b^T \nabla \eta\right)\\
&&-(I-\Delta)^{-1}\text{div}\left(\frac{1}{2}(\nabla(u+v):\nabla a+\nabla(\gamma+\eta)\cdot\nabla b+(\gamma+\eta)b) I \right)\\
&\triangleq & I_1+I_2+I_3+I_4,
\end{eqnarray*}
and
\begin{eqnarray*}
R_2(t,x)&\triangleq &
-a\cdot\nabla \eta-(I-\Delta)^{-1}\left(b(\text{div} v)+\gamma(\text{div} a)\right)\\
&&-(I-\Delta)^{-1}\text{div}\left(\nabla b \nabla v+(\nabla b)\cdot \nabla v-\nabla b(\text{div} v)\right)\\
&&-(I-\Delta)^{-1}\text{div}\left(\nabla \gamma \nabla a+(\nabla \gamma)\cdot \nabla a-\nabla \gamma(\text{div} a)\right).
\end{eqnarray*}

We first claim that for all $s>\max (1+\frac{d}{p},\frac 3 2)$ and $t\in [0,T]$,
\begin{eqnarray}\label{3claim}
\quad\quad ||R_1(t)||_{B^{s-1}_{p,r}}+||R_2(t)||_{B^{s-1}_{p,r}}
\leq C(||a(t)||_{B^{s-1}_{p,r}}+||b(t)||_{B^{s-1}_{p,r}}) A(t;s),
\end{eqnarray}
where $A(t;s)\triangleq ||u(t)||_{B^s_{p,r}}+||\gamma(t)||_{B^s_{p,r}}+||v(t)||_{B^s_{p,r}}+||\eta(t)||_{B^s_{p,r}}$.\\
Indeed, for $s>1+\frac{d}{p}$, $B^{s-1}_{p,r}(\mathbb{R}^d)$ is an algebra, one has
\begin{eqnarray*}
||-a\cdot\nabla v||_{B^{s-1}_{p,r}}
\leq C||a||_{B^{s-1}_{p,r}}||\nabla v||_{B^{s-1}_{p,r}}
\leq C||a||_{B^{s-1}_{p,r}}||v||_{B^s_{p,r}}.
\end{eqnarray*}
\textbf{Case 1: $\max (1+\frac{d}{p},\frac 3 2)<s\leq 2+\frac{d}{p}.$}\\
Since both $-(I-\Delta)^{-1}\text{div}$ and $-(I-\Delta)^{-1}$ belong to $Op(S^{-1})$, it follows from
Lemma \ref{2-multiplier} and (\ref{2morse2}) that
\begin{eqnarray*}
&&||-(I-\Delta)^{-1}\left(a(\text{div} v)+u(\text{div} a)+a\cdot \nabla v^T+u\cdot \nabla a^T\right)||_{B^{s-1}_{p,r}}\\ \nonumber
&\leq& C ||a(\text{div} v)+u(\text{div} a)+a\cdot \nabla v^T+u\cdot \nabla a^T||_{B^{s-2}_{p,r}}\\ \nonumber
&\leq& C (||a||_{B^{s-1}_{p,r}} ||\nabla v||_{B^{s-2}_{p,r}}+||u||_{B^{s-1}_{p,r}} ||\nabla a||_{B^{s-2}_{p,r}}) \\ \nonumber
&\leq& C||a||_{B^{s-1}_{p,r}}(||u||_{B^s_{p,r}}+||v||_{B^s_{p,r}}),
\end{eqnarray*}
which yields
\begin{eqnarray*}
||I_1||_{B^{s-1}_{p,r}}
\leq C||a||_{B^{s-1}_{p,r}}(||u||_{B^s_{p,r}}+||v||_{B^s_{p,r}}).
\end{eqnarray*}

\begin{eqnarray*}
||I_3||_{B^{s-1}_{p,r}}
&\leq& C ||-\nabla a (\text{div} v)-\nabla u (\text{div} a)-\nabla \gamma^T \nabla b-\nabla b^T \nabla \eta||_{B^{s-2}_{p,r}}\\ \nonumber
&\leq& C ||\nabla a||_{B^{s-2}_{p,r}} (||\nabla u||_{B^{s-1}_{p,r}}+||\text{div} v||_{B^{s-1}_{p,r}})\\ \nonumber
&&+ C ||\nabla b||_{B^{s-2}_{p,r}} (||\nabla \gamma||_{B^{s-1}_{p,r}}+||\nabla \eta||_{B^{s-1}_{p,r}})\\ \nonumber
&\leq& C||a||_{B^{s-1}_{p,r}}(||u||_{B^s_{p,r}}+||v||_{B^s_{p,r}})
+ C||b||_{B^{s-1}_{p,r}}(||\gamma||_{B^s_{p,r}}+||\eta||_{B^s_{p,r}}).
\end{eqnarray*}
Likewise,
\begin{eqnarray*}
||I_2+I_4||_{B^{s-1}_{p,r}}
\leq C||a||_{B^{s-1}_{p,r}}(||u||_{B^s_{p,r}}+||v||_{B^s_{p,r}})
+ C||b||_{B^{s-1}_{p,r}}(||\gamma||_{B^s_{p,r}}+||\eta||_{B^s_{p,r}}),
\end{eqnarray*}
which implies
\begin{eqnarray}\label{3R-1}
||R_1(t)||_{B^{s-1}_{p,r}}
&\leq& C||a||_{B^{s-1}_{p,r}}(||u||_{B^s_{p,r}}+||v||_{B^s_{p,r}})\\ \nonumber
&&+ C||b||_{B^{s-1}_{p,r}}(||\gamma||_{B^s_{p,r}}+||\eta||_{B^s_{p,r}}).
\end{eqnarray}
Similarly,
\begin{eqnarray}\label{3R-2}
\quad\quad\quad ||R_2(t)||_{B^{s-1}_{p,r}}
\leq C||a||_{B^{s-1}_{p,r}}(||\gamma||_{B^s_{p,r}}+||\eta||_{B^s_{p,r}})+ C||b||_{B^{s-1}_{p,r}}||v||_{B^s_{p,r}}.
\end{eqnarray}
\textbf{Case 2: $s>2+\frac{d}{p}.$}
Notice that $B^{s-2}_{p,r}(\mathbb{R}^d)$ is an algebra, which ensures that (\ref{3R-1}) and (\ref{3R-2}) still hold true.
Thus, we have proven the above claim (\ref{3claim}).

On the other hand, for $s>1+\frac{d}{p}$, we have
$$||\nabla u||_{B^{d\over p}_{p,r}\cap L^\infty}\leq C ||u||_{B^{s}_{p,r}}\,\ \,\ \text{and}\,\ \,\
||\nabla u||_{B^{s-2}_{p,r}}\leq C ||u||_{B^{s}_{p,r}}.$$
Applying Lemma \ref{2transport-estimates} (i) to system (\ref{3diffence-eqs}) yields, for the case (i),
\begin{eqnarray*}
||a(t)||_{B^{s-1}_{p,r}}
\leq ||a_0||_{B^{s-1}_{p,r}}+\int_0^t ||R_1(\tau)||_{B^{s-1}_{p,r}}d\tau
+C \int_0^t ||u(\tau)||_{B^s_{p,r}}||a(\tau)||_{B^{s-1}_{p,r}}d\tau
\end{eqnarray*}
and
\begin{eqnarray*}
||b(t)||_{B^{s-1}_{p,r}}
\leq ||b_0||_{B^{s-1}_{p,r}}+\int_0^t ||R_2(\tau)||_{B^{s-1}_{p,r}}d\tau
+C \int_0^t ||u(\tau)||_{B^s_{p,r}}||b(\tau)||_{B^{s-1}_{p,r}}d\tau,
\end{eqnarray*}
which together with (\ref{3claim}) leads to
\begin{eqnarray*}
&& ||a(t)||_{B^{s-1}_{p,r}}+||b(t)||_{B^{s-1}_{p,r}}\\
&\leq& ||a_0||_{B^{s-1}_{p,r}}+||b_0||_{B^{s-1}_{p,r}}
+C \int_0^t (||a(\tau)||_{B^{s-1}_{p,r}}+||b(\tau)||_{B^{s-1}_{p,r}}) A(\tau;s) d\tau.
\end{eqnarray*}
Taking advantage of the Gronwall inequality, one reaches (\ref{3-lemma1-i}).

For the critical case (ii), we here use the interpolation method to handle it.
Indeed, if we choose
$s_1\in(\max(1+\frac{d}{p},\frac 3 2)-1,s-1)$, $s_2\in (s-1,s)$ and $\theta=\frac{s_2-(s-1)}{s_2-s_1} \in(0,1)$,
then $s-1=\theta s_1+(1-\theta) s_2$.
According to (\ref{2interpolation}) and (\ref{3-lemma1-i}), one deduces that
\begin{eqnarray*}
&& ||a(t)||_{B^{s-1}_{p,r}}+||b(t)||_{B^{s-1}_{p,r}}\\
&\leq& (||a(t)||_{B^{s_1}_{p,r}}+||b(t)||_{B^{s_1}_{p,r}})^{\theta}  (||a(t)||_{B^{s_2}_{p,r}}+||b(t)||_{B^{s_2}_{p,r}})^{1-\theta}\\ \nonumber
&\leq& C M^\theta(t;s_1)A^{1-\theta}(t;s_2)\\ \nonumber
&\leq& C M^\theta(t;s-1)A^{1-\theta}(t;s).
\end{eqnarray*}
Therefore, we complete our proof of Lemma 3.1.
\end{proof}

Next, we construct the approximation solutions to system (\ref{1mchnew}) as follows.
\begin{lemma3}
Let $d\in \mathbb{N_+}$, $1\leq p,r\leq \infty$ and $s>1+\frac{d}{p}$
(or $s=1+\frac{d}{p}$ with $r=1$ and $1\leq p<\infty$).
Assume that $(u_0,\gamma_0)\in B^s_{p,r}(\mathbb{R}^d)\times B^s_{p,r}(\mathbb{R}^d)$ and $u^0=\gamma^0\equiv 0$.
Then\\
(i) there exists a sequence of smooth functions
$(u^n,\gamma^n)_{n\in \mathbb{N}}\in C(\mathbb{R}^{+}; B^{\infty}_{p,r}(\mathbb{R}^d)\times B^{\infty}_{p,r}(\mathbb{R}^d))$
which solves the following linear transport equations by induction with respect to $n$:
\[(TE_n)\left\{\begin{array}{l}
(\partial_t+u^n\cdot \nabla) u^{n+1}=F_1(u^n,\gamma^n)\triangleq F^n_1(t,x),\\
(\partial_t+u^n\cdot \nabla) \gamma^{n+1}=F_2(u^n,\gamma^n)\triangleq F^n_2(t,x),\\
u^{n+1}|_{t=0}\triangleq u^{n+1}_0(x)=S_{n+1} u_0(x),\\
\gamma^{n+1}|_{t=0}\triangleq \gamma^{n+1}_0(x)=S_{n+1} \gamma_0(x),
\end{array}\right.\]
where $F_1(u^n,\gamma^n), F_2(u^n,\gamma^n)$ are defined by (\ref{1-F1}) and (\ref{1-F2}),
and $S_{n+1}\triangleq {\sum\limits_{q=-1}^{n}}\Delta_q$ is the low frequency cut-off operator.\\
(ii) there exists a time $T>0$ such that the solution $(u^n,\gamma^n)_{n\in \mathbb{N}}$
is uniformly bounded in $E^s_{p,r}(T)\times E^s_{p,r}(T)$.\\
(iii) if we further suppose that $s>\max(1+\frac{d}{p},\frac{3}{2})$,
then $(u^n,\gamma^n)_{n\in \mathbb{N}}$ is  a Cauchy sequence in $C([0,T]; B^{s-1}_{p,r}(\mathbb{R}^d)\times B^{s-1}_{p,r}(\mathbb{R}^d))$
and thus converges to a limit $(u,\gamma)\in C([0,T]; B^{s-1}_{p,r}(\mathbb{R}^d)\times B^{s-1}_{p,r}(\mathbb{R}^d))$.
\end{lemma3}

\begin{proof}
(i) Thanks to all the data $S_{n+1}u_0 \in B^{\infty}_{p,r}(\mathbb{R}^d)$, by induction with respect to the index $n$
and applying the existence and uniqueness theory of transport equations \cite{BCD} to $(TE_n)$, one can easily get the desired result.

(ii) Applying Lemma \ref{2transport-estimates} (i) to $(TE_n)$, one gets
\begin{eqnarray*}
||u^{n+1}(t)||_{B^{s}_{p,r}}
&\leq& ||S_{n+1}u_0||_{B^{s}_{p,r}}
+C \int_0^t ||\nabla u^n(\tau)||_{B^{s-1}_{p,r}}||u^{n+1}(\tau)||_{B^{s}_{p,r}}d\tau\\
&&+\int_0^t ||F^n_1(\tau)||_{B^{s}_{p,r}}d\tau,
\end{eqnarray*}
and
\begin{eqnarray*}
||\gamma^{n+1}(t)||_{B^{s}_{p,r}}
&\leq& ||S_{n+1}\gamma_0||_{B^{s}_{p,r}}
+C \int_0^t ||\nabla u^n(\tau)||_{B^{s-1}_{p,r}}||\gamma^{n+1}(\tau)||_{B^{s}_{p,r}}d\tau\\
&&+\int_0^t ||F^n_2(\tau)||_{B^{s}_{p,r}}d\tau.
\end{eqnarray*}
On the other hand, noting that $B^{s-1}_{p,r}(\mathbb{R}^d)$ is an algebra and simulating the proof of (\ref{3claim}), one obtains
\begin{eqnarray*}
||F^n_1(t)||_{B^{s}_{p,r}}+||F^n_2(t)||_{B^{s}_{p,r}}
\leq C (||u^n(t)||_{B^{s}_{p,r}}+||\gamma^n(t)||_{B^{s}_{p,r}})^2.
\end{eqnarray*}
Set
\begin{eqnarray*}
\Gamma^n(t)\triangleq ||u^n(t,\cdot)||_{B^{s}_{p,r}}+||\gamma^n(t,\cdot)||_{B^{s}_{p,r}}.
\end{eqnarray*}
Then the Gronwall inequality gives
\begin{eqnarray}\label{3Gamma-estimates}
\quad\quad \Gamma^{n+1}(t)
\leq C e^{C U^n(t)}
\left(||u_0||_{B^{s}_{p,r}}+||\gamma_0||_{B^{s}_{p,r}}+\int_0^t e^{-C U^n(\tau)}(\Gamma^n(\tau))^2d\tau\right)
\end{eqnarray}
with $U^n(t)\triangleq \int_0^t ||u^n(\tau)||_{B^s_{p,r}}d\tau$.\\
Choose $0<\, T\, < \frac{1}{2C^2 (||u_0||_{B^{s}_{p,r}}+||\gamma_0||_{B^{s}_{p,r}})}$ and suppose that
\begin{eqnarray}\label{3Gamma-bdd}
\Gamma^n(t) \leq \frac{C(||u_0||_{B^{s}_{p,r}}+||\gamma_0||_{B^{s}_{p,r}})}{1-2C^2 (||u_0||_{B^{s}_{p,r}}+||\gamma_0||_{B^{s}_{p,r}}) t},\quad \forall t\in[0,T].
\end{eqnarray}
Noting that $e^{C (U^n(t)-U^n(\tau))}\leq
\sqrt{\frac{1-2C^2 (||u_0||_{B^{s}_{p,r}}+||\gamma_0||_{B^{s}_{p,r}})\tau}{1-2C^2 (||u_0||_{B^{s}_{p,r}}+||\gamma_0||_{B^{s}_{p,r}})t}}$
and substituting (\ref{3Gamma-bdd}) into (\ref{3Gamma-estimates}) yields
\begin{eqnarray*}
\Gamma^{n+1}(t)
&\leq&\frac{C(||u_0||_{B^{s}_{p,r}}+||\gamma_0||_{B^{s}_{p,r}})}{\sqrt{1-2C^2(||u_0||_{B^{s}_{p,r}}+||\gamma_0||_{B^{s}_{p,r}})t}}
+\frac{C}{\sqrt{1-2C^2 (||u_0||_{B^{s}_{p,r}}+||\gamma_0||_{B^{s}_{p,r}})t}}\\
&&\times \int_0^t \frac{C^2 (||u_0||_{B^{s}_{p,r}}+||\gamma_0||_{B^{s}_{p,r}})^2}{\left(1-2C^2 (||u_0||_{B^{s}_{p,r}}+||\gamma_0||_{B^{s}_{p,r}})\tau\right)^\frac{3}{2}} d\tau\\
\nonumber&=&\frac{C(||u_0||_{B^{s}_{p,r}}+||\gamma_0||_{B^{s}_{p,r}})}{\sqrt{1-2C^2 (||u_0||_{B^{s}_{p,r}}+||\gamma_0||_{B^{s}_{p,r}})t}}
+\frac{C (||u_0||_{B^{s}_{p,r}}+||\gamma_0||_{B^{s}_{p,r}})}{\sqrt{1-2C^2 (||u_0||_{B^{s}_{p,r}}+||\gamma_0||_{B^{s}_{p,r}})t}}\\
&&\times\left(\frac{1}{\sqrt{1-2C^2 (||u_0||_{B^{s}_{p,r}}+||\gamma_0||_{B^{s}_{p,r}})t}}-1\right)\\
\nonumber&\leq&\frac{C(||u_0||_{B^{s}_{p,r}}+||\gamma_0||_{B^{s}_{p,r}})}{1-2C^2 (||u_0||_{B^{s}_{p,r}}+||\gamma_0||_{B^{s}_{p,r}})t},
\end{eqnarray*}
which implies that $(u^n,\gamma^n)_{n\in \mathbb{N}}$ is uniformly bounded in $C([0,T]; B^s_{p,r}\times B^s_{p,r})$.
By using system $(TE_n)$ and the similar proof of (\ref{3claim}), one can readily deduce that
$(\partial_t u^{n+1},\partial_t \gamma^{n+1})_{n\in \mathbb{N}}$ is uniformly bounded in $C([0,T]; B^{s-1}_{p,r}\times B^{s-1}_{p,r})$.
Hence, we have proven (ii).

(iii) For all $m,n\in \mathbb{N}$, by system $(TE_n)$ again, we have
\[\left\{\begin{array}{l}
(\partial_t+u^{n+m}\cdot \nabla) (u^{n+m+1}-u^{n+1})
=F^{n+m}_1(t,x)-F^n_1(t,x)+(u^{n}-u^{n+m})\cdot \nabla u^{n+1},\\
(\partial_t+u^{n+m}\cdot \nabla) (\gamma^{n+m+1}-\gamma^{n+1})
=F^{n+m}_2(t,x)-F^n_2(t,x)+(u^{n}-u^{n+m})\cdot \nabla \gamma^{n+1},
\end{array}\right.\]
where $F^{n+m}_i(t,x)\triangleq F_i(u^{n+m},\gamma^{n+m})$ ($i=1,2$) are defined by (\ref{1-F1}) and (\ref{1-F2}).\\
Similar to the proof of (\ref{3-lemma1-i}), for $s>\max (1+\frac{d}{p},\frac 3 2)$ and $s\neq 2+\frac{d}{p}$,
one gets
\begin{eqnarray*}
\Lambda^{n+m+1}_{n+1}(t)
&\leq& C e^{C U^{n+m}(t)}(||u_0^{n+m+1}-u_0^{n+1}||_{B^{s-1}_{p,r}}+||\gamma_0^{n+m+1}-\gamma_0^{n+1}||_{B^{s-1}_{p,r}}\\
&&+\int_0^t e^{-C U^{n+m}(\tau)} \Lambda^{n+m}_{n}(\tau) (\Gamma^{n}(\tau)+\Gamma^{n+1}(\tau)+\Gamma^{n+m}(\tau))d\tau),
\end{eqnarray*}
where $\Lambda^{n+m}_{n}(t)\triangleq ||(u^{n+m}-u^{n})(t)||_{B^{s-1}_{p,r}}+||(\gamma^{n+m}-\gamma^{n})(t)||_{B^{s-1}_{p,r}}$.\\
Note that
\begin{eqnarray*}
||u_0^{n+m+1}-u_0^{n+1}||_{B^{s-1}_{p,r}}+||\gamma_0^{n+m+1}-\gamma_0^{n+1}||_{B^{s-1}_{p,r}}
\leq C 2^{-n}(||u_0||_{B^s_{p,r}}+||\gamma_0||_{B^s_{p,r}}).
\end{eqnarray*}
Then according to (ii), one can find a constant $C_{T}>0$, independent of $n,m$, such that for all $t\in[0,T]$,
\begin{eqnarray*}
\Lambda^{n+m+1}_{n+1}(t)
\leq C_{T}\left(2^{-n}+\int_0^t  \Lambda^{n+m}_{n}(\tau) d\tau\right).
\end{eqnarray*}
Arguing by induction with respect to the index $n$, we have
\begin{eqnarray*}
\Lambda^{n+m+1}_{n+1}(t)
&\leq& C_{T}\left(2^{-n}\sum\limits_{k=0}^n \frac{(2T C_T)^k}{k!}+C^{n+1}_T \int_0^t  \frac{(t-\tau)^n}{n!}d\tau\right)\\
&\leq& \left(C_{T}\sum\limits_{k=0}^n \frac{(2T C_T)^k}{k!}\right)2^{-n}+C_T \frac{(T C_T)^{n+1}}{(n+1)!},
\end{eqnarray*}
which implies the desired result.

While for the critical point $s=2+\frac{d}{p}$,
we can apply the similar interpolation argument used in the proof of (\ref{3-lemma1-ii}) to show that
$(u^n,\gamma^n)_{n\in \mathbb{N}}$ is also a Cauchy sequence in $C([0,T]; B^{s-1}_{p,r}\times B^{s-1}_{p,r})$.
Therefore, we have completed the proof of Lemma 3.2.
\end{proof}

\textbf{Proof of Theorem 1.1.}
We first claim that the obtained limit $(u,\gamma)$ in Lemma 3.2 (iii) belongs to $E^s_{p,r}(T)\times E^s_{p,r}(T)$ and solves system (\ref{1mchnew}).
In fact, according to Lemma 3.2 (ii) and the Fatou lemma, we have $(u,\gamma)\in L^{\infty}(0,T; B^s_{p,r}\times B^s_{p,r})$.

By Lemma 3.2 (iii) again, then an interpolation argument gives
$$(u^n,\gamma^n) \rightarrow (u,\gamma)\,\ \text{in}\,\ C([0,T]; B^{s'}_{p,r}\times B^{s'}_{p,r}),\,\ \text{as} \,\ n\to \infty, \,\ \forall\, s'<s.$$
Then taking limit in $(TE_n)$, one can see that $(u,\gamma)$ solves system (\ref{1mchnew}) in the sense of $C([0,T]; B^{s'-1}_{p,r}\times B^{s'-1}_{p,r})$ for all $s'<s$.

In view of $(u,\gamma)\in L^{\infty}(0,T; B^s_{p,r}\times B^s_{p,r})$ and similar to the proof of (\ref{3claim}),
$F_1(u,\gamma)$ and $F_2(u,\gamma)$ in system (\ref{1mchnew}) also belong to $L^{\infty}(0,T; B^s_{p,r}\times B^s_{p,r})$.
While $u\cdot\nabla u$ and $u\cdot\nabla \gamma$ belong to $L^{\infty}(0,T; B^{s-1}_{p,r}\times B^{s-1}_{p,r})$,
thus $(\partial_t u,\partial_t \gamma)\in L^{\infty}(0,T; B^{s-1}_{p,r}\times B^{s-1}_{p,r})$.
Furthermore, if $r<\infty$, then Lemma \ref{2transport-estimates} (ii) ensures $(u,\gamma)\in C([0,T]; B^s_{p,r}\times B^s_{p,r})$.
Making use of system (\ref{1mch}) again, one obtains $(\partial_t u,\partial_t \gamma)\in C([0,T]; B^{s-1}_{p,r}\times B^{s-1}_{p,r})$.
Hence, $(u,\gamma)\in E^s_{p,r}(T)\times E^s_{p,r}(T)$.

On the other hand, the continuity with respect to the initial data in
$$ C([0,T]; B^{s'}_{p,r}\times B^{s'}_{p,r})\cap C^1{([0,T]; B^{s'-1}_{p,r}\times B^{s'-1}_{p,r})} \quad (\forall\, s'<s)$$
can be easily proved by Lemma 3.1 and an interpolation argument.
While the continuity up to $s'=s$ in the case of $r<\infty$
can be obtained through the use of a sequence of viscosity approximation solutions
$(u_{\varepsilon},\gamma_{\varepsilon})_{\varepsilon>0}$ for system (\ref{1mchnew}) which converges uniformly in
$C([0,T]; B^{s}_{p,r}\times B^{s}_{p,r})\cap C^1{([0,T]; B^{s-1}_{p,r}\times B^{s-1}_{p,r})}$.
Therefore, we have proven Theorem 1.1.
\quad \quad \quad \quad \quad \quad \quad \quad \quad \quad \quad \quad \quad\quad \quad \quad \quad \quad $\square$

\section{Local well-posedness in critical Besov spaces}
\newtheorem {remark4}{Remark}[section]
\newtheorem{theorem4}{Theorem}[section]
\newtheorem{lemma4}{Lemma}[section]
\newtheorem{corollary4}{Corollary}[section]
\newtheorem{proposition4}{Proposition}[section]

In this section, we shall establish the local well-posedness of system (\ref{1mchnew}) in critical Besov spaces.
In order to prove Theorem 1.2, let us first give the existence of solutions as follows.
\begin{lemma4}
Let $d\in \mathbb{N_+}$ and  $1\leq p< \infty$.
Suppose that $(u_0,\gamma_0)\in B^{1+\frac{d}{p}}_{p,1}(\mathbb{R}^d)\times B^{1+\frac{d}{p}}_{p,1}(\mathbb{R}^d)$.
Then there exists a time $T>0$ such that system (\ref{1mchnew}) has a solution
$(u,\gamma)\in E^{1+\frac{d}{p}}_{p,1}(T)\times E^{1+\frac{d}{p}}_{p,1}(T)$.
Moreover, for some fixed $\delta>0$,  there exists a constant
$M=M(\delta,||(u_0,\gamma_0)||_{B^{1+\frac{d}{p}}_{p,1}(\mathbb{R}^d)\times B^{1+\frac{d}{p}}_{p,1}(\mathbb{R}^d)})>0$ such that for all
$(v_0,\eta_0)\in B^{1+\frac{d}{p}}_{p,1}(\mathbb{R}^d)\times B^{1+\frac{d}{p}}_{p,1}(\mathbb{R}^d)$ with
$||(v_0-u_0, \eta_0-\gamma_0)||_{B^{1+\frac{d}{p}}_{p,1}\times B^{1+\frac{d}{p}}_{p,1}} \leq \delta $,
the system (\ref{1mchnew}) has a solution $(v,\eta)\in E^{1+\frac{d}{p}}_{p,1}(T)\times E^{1+\frac{d}{p}}_{p,1}(T)$
satisfying $||(v,\eta)||_{L^\infty (0,T; B^{1+\frac{d}{p}}_{p,1}(\mathbb{R}^d)\times B^{1+\frac{d}{p}}_{p,1}(\mathbb{R}^d))}\leq M $.
\end{lemma4}

\begin{proof}
Thanks to Lemma 3.2 (ii), the smooth approximation solution $(u^n,\gamma^n)$ to $(TE_n)$ is uniformly bounded in
$E^{1+\frac{d}{p}}_{p,1}(T)\times E^{1+\frac{d}{p}}_{p,1}(T)$.
Then the Arzela-Ascoli theorem and a standard diagonal process ensures that, up to an extraction,
$(u^n,\gamma^n)$ converges to a limit $(u,\gamma)$ in $C([0,T]; (B^{\frac{d}{p}}_{p,1})_{loc}\times (B^{\frac{d}{p}}_{p,1})_{loc})$.
Besides, by Lemma 3.2 (ii) and the Fatou lemma,
we get $(u,\gamma)\in L^\infty(0,T; B^{1+\frac{d}{p}}_{p,1}\times B^{1+\frac{d}{p}}_{p,1})$.
This together with an interpolation argument leads to
$(u^n,\gamma^n)$ converges to $(u,\gamma)$ in
$C([0,T]; (B^{s}_{p,1})_{loc}\times (B^{s}_{p,1})_{loc})$ for any $s<1+\frac{d}{p}$.
Then taking limit in $(TE_n)$, one deduces that $(u,\gamma)$ is indeed a solution to system (\ref{1mchnew}).

On the other hand, since $(u,\gamma)\in L^\infty(0,T; B^{1+\frac{d}{p}}_{p,1}\times B^{1+\frac{d}{p}}_{p,1})$,
it then follows from Lemma \ref{2transport-estimates} (ii) that
$(u,\gamma)\in C([0,T]; B^{1+\frac{d}{p}}_{p,1}\times B^{1+\frac{d}{p}}_{p,1})$.
By using system (\ref{1mchnew}) again, one can readily infer that
$(\partial_t u, \partial_t \gamma)\in C([0,T]; B^{\frac{d}{p}}_{p,1}\times B^{\frac{d}{p}}_{p,1})$.
Thus, $(u,\gamma)\in E^{1+\frac{d}{p}}_{p,1}(T)\times E^{1+\frac{d}{p}}_{p,1}(T)$.

Next, by the assumption,
$||v_0||_{B^{1+\frac{d}{p}}_{p,1}}+||\eta_0||_{B^{1+\frac{d}{p}}_{p,1}}
\leq ||u_0||_{B^{1+\frac{d}{p}}_{p,1}}+||\gamma_0||_{B^{1+\frac{d}{p}}_{p,1}} +\delta$.
According to (\ref{3Gamma-bdd}), by choosing
$T\triangleq \frac{1}{4C^2 \left(||u_0||_{B^{1+\frac{d}{p}}_{p,1}}+||\gamma_0||_{B^{1+\frac{d}{p}}_{p,1}}+ \delta\right)}>0$
and
$M\triangleq 2(||u_0||_{B^{1+\frac{d}{p}}_{p,1}}+ ||\gamma_0||_{B^{1+\frac{d}{p}}_{p,1}}+\delta)>0$,
one can easily complete the proof of Lemma 4.1.
\end{proof}

With regard to uniqueness and continuous dependence of the solution,
one cannot directly follow the similar lines of the proof in Theorem 1.1 to verify Theorem 1.2,
since Lemma 3.1 fails in the case of $s=1+\frac{d}{p}$ and $r=1$, while it is the cornerstone in the proof of Theorem 1.1.
To overcome this difficulty, we first notice the following interpolation inequality
\begin{eqnarray} \label{4interpolation}
||f||_{B^{\frac{d}{p}}_{p,1}}
\leq ||f||^{\theta}_{B^{1+\frac{d}{p}}_{p,1}} ||f||^{1-\theta}_{B^{d/p-\theta/ (1-\theta)}_{p,1}}
\leq C ||f||^{\theta}_{B^{1+\frac{d}{p}}_{p,1}} ||f||^{1-\theta}_{B^{\frac{d}{p}}_{p,\infty}}
\end{eqnarray}
with $\theta\in(0,1)$.
Furthermore, we can establish  a priori estimates for the solution in
$L^\infty(0,T; B^{\frac{d}{p}}_{p,\infty}(\mathbb{R}^d)\times B^{\frac{d}{p}}_{p,\infty}(\mathbb{R}^d))$ as $1\leq p\leq 2d$ (see Lemma 4.2 below),
which together with the uniform bounds for the solution in $C([0,T]; B^{1+\frac{d}{p}}_{p,1}(\mathbb{R}^d)\times B^{1+\frac{d}{p}}_{p,1}(\mathbb{R}^d))$
yields the continuity in $C([0,T]; B^{\frac{d}{p}}_{p,1}(\mathbb{R}^d)\times B^{\frac{d}{p}}_{p,1}(\mathbb{R}^d))$.

\begin{lemma4}
Let $d\in \mathbb{N_+}$ and  $1\leq p\leq 2d$. Assume that
$(u,\gamma), (v,\eta) \in L^{\infty}(0,T; (B^{1+\frac{d}{p}}_{p,\infty}\cap Lip)\times (B^{1+\frac{d}{p}}_{p,\infty}\cap Lip))
\cap C([0,T]; B^{\frac{d}{p}}_{p,\infty}\times B^{\frac{d}{p}}_{p,\infty})$ are two solutions to system (\ref{1mchnew})
with the initial data $(u_0,\gamma_0), (v_0,\eta_0)\in (B^{1+\frac{d}{p}}_{p,\infty}\cap Lip)\times (B^{1+\frac{d}{p}}_{p,\infty}\cap Lip)$, respectively.
Set $(a,b)\triangleq (v-u,\eta-\gamma)$ and $(a_0,b_0)\triangleq (v_0-u_0,\eta_0-\gamma_0)$.
Let
\begin{eqnarray*}
H(t) \triangleq
e^{-C\int_0^t ||\nabla u(\tau)||_{B^{\frac{d}{p}}_{p,\infty}\cap L^{\infty}}d\tau}
\left(||a(t)||_{B^{\frac{d}{p}}_{p,\infty}(\mathbb{R}^d)}+||b(t)||_{B^{\frac{d}{p}}_{p,\infty}(\mathbb{R}^d)}\right)
\end{eqnarray*}
and
\begin{eqnarray*}
G(t) \triangleq
||u(t)||_{B^{1+\frac{d}{p}}_{p,\infty}\cap Lip}+||\gamma(t)||_{B^{1+\frac{d}{p}}_{p,\infty}\cap Lip}
+||v(t)||_{B^{1+\frac{d}{p}}_{p,\infty}\cap Lip}+||\eta(t)||_{B^{1+\frac{d}{p}}_{p,\infty}\cap Lip}.
\end{eqnarray*}
If there exists a constant $C>0$ such that for any $T^\star\leq T$,
\begin{eqnarray}\label{4lemma2-assumption}
\sup\limits_{t\in[0,T^\star]} H(t)\leq 1,
\end{eqnarray}
then for all $t\in [0,T^{\star}]$, we have
\begin{eqnarray}\label{4lemma2-estimates}
&&\frac{||a(t)||_{B^{\frac{d}{p}}_{p,\infty}(\mathbb{R}^d)}+||b(t)||_{B^{\frac{d}{p}}_{p,\infty}(\mathbb{R}^d)}}{e}\\ \nonumber
&\leq& e^{C\int_0^t ||\nabla u(\tau)||_{B^{\frac{d}{p}}_{p,\infty}\cap L^{\infty}}d\tau}
\left(\frac{||a_0||_{B^{\frac{d}{p}}_{p,\infty}(\mathbb{R}^d)}+||b_0||_{B^{\frac{d}{p}}_{p,\infty}(\mathbb{R}^d)}}{e}\right)^{\exp(-C\int_0^t L(G(\tau))d\tau)}
\end{eqnarray}
with $L(x)\triangleq x\ln(e+Kx)$ and
$K\triangleq \sup\limits_{t\in[0,T]} \left(\frac{||a(t)||_{B^{\frac{d}{p}}_{p,\infty}}}{||b(t)||_{B^{\frac{d}{p}}_{p,\infty}}}
+\frac{||b(t)||_{B^{\frac{d}{p}}_{p,\infty}}}{||a(t)||_{B^{\frac{d}{p}}_{p,\infty}}}+2\right)$.\\
In particular, (\ref{4lemma2-estimates}) is also true on $[0,T]$ provided that
\begin{eqnarray}\label{4lemma2-initial-assumption}
||a_0||_{B^{\frac{d}{p}}_{p,\infty}(\mathbb{R}^d)}+||b_0||_{B^{\frac{d}{p}}_{p,\infty}(\mathbb{R}^d)}
\leq e^{1-\exp(C\int_0^T L(G(t))d t)}.
\end{eqnarray}
\end{lemma4}

\begin{proof}
We first declare that for any $1\leq p\leq 2d$, $t\in [0,T]$, we have
\begin{eqnarray}\label{4lemma2-proof-RHS-estimates}
&&||R_1(t,\cdot)||_{B^{\frac{d}{p}}_{p,\infty}(\mathbb{R}^d)}+||R_2(t,\cdot)||_{B^{\frac{d}{p}}_{p,\infty}(\mathbb{R}^d)}\\ \nonumber
&\leq& C \left(||a(t)||_{B^{\frac{d}{p}}_{p,1}(\mathbb{R}^d)}+||b(t)||_{B^{\frac{d}{p}}_{p,1}(\mathbb{R}^d)}\right) G(t),
\end{eqnarray}
where $R_1(t,x)$ and $R_2(t,x)$ are defined in (\ref{3diffence-eqs}).\\
Indeed, since both $B^{\frac{d}{p}}_{p,\infty}(\mathbb{R}^d)\cap L^\infty(\mathbb{R}^d)$ and
$B^{\frac{d}{p}}_{p,1}(\mathbb{R}^d)$ are algebras, then
\begin{eqnarray}\label{4lemma2-proof-RHS-estimates-1}
&& ||-a\cdot\nabla v||_{B^{\frac{d}{p}}_{p,\infty}}+||-a\cdot\nabla \eta||_{B^{\frac{d}{p}}_{p,\infty}}\\ \nonumber
&\leq& C||a||_{B^{\frac{d}{p}}_{p,\infty}\cap L^\infty}
\left(||\nabla v||_{B^{\frac{d}{p}}_{p,\infty}\cap L^\infty}+||\nabla \eta||_{B^{\frac{d}{p}}_{p,\infty}\cap L^\infty}\right)\\ \nonumber
&\leq& C||a||_{B^{\frac{d}{p}}_{p,1}} \left(||v||_{B^{1+\frac{d}{p}}_{p,\infty}\cap Lip}+||\eta||_{B^{1+\frac{d}{p}}_{p,\infty}\cap Lip}\right).
\end{eqnarray}
Noting that $-(I-\Delta)^{-1}\text{div}, -(I-\Delta)^{-1}\in Op(S^{-1})$,
then applying Lemma \ref{2-multiplier} and (\ref{2morsecritical}), one can easily get for $1\leq p\leq 2d$,
\begin{eqnarray*}
&& ||R_1(t,\cdot)+a\cdot\nabla v||_{B^{\frac{d}{p}}_{p,\infty}}+||R_2(t,\cdot)+a\cdot\nabla \eta||_{B^{\frac{d}{p}}_{p,\infty}}\\
&\leq& C \left(||a(t)||_{B^{\frac{d}{p}}_{p,1}}+||b(t)||_{B^{\frac{d}{p}}_{p,1}}\right) G(t),
\end{eqnarray*}
which along with (\ref{4lemma2-proof-RHS-estimates-1}) leads to (\ref{4lemma2-proof-RHS-estimates}).

Applying Lemma \ref{2transport-estimates} (i) to (\ref{3diffence-eqs}), one infers
\begin{eqnarray*}
H(t) \leq H(0)+ \int_0^t e^{-C\int_0^{\tau} ||\nabla u(\xi)||_{B^{\frac{d}{p}}_{p,\infty}\cap L^{\infty}}d \xi}
\left(||R_1(\tau)||_{B^{\frac{d}{p}}_{p,\infty}}+||R_2(\tau)||_{B^{\frac{d}{p}}_{p,\infty}}\right)d\tau,
\end{eqnarray*}
which together with (\ref{4lemma2-proof-RHS-estimates}) yields
\begin{eqnarray} \label{4lemma2-proof-H-estimates-1}
H(t) &\leq& H(0)+C \int_0^t G(\tau)e^{-C\int_0^{\tau} ||\nabla u(\xi)||_{B^{\frac{d}{p}}_{p,\infty}\cap L^{\infty}}d \xi}\\ \nonumber
&&\times \left(||a(\tau)||_{B^{\frac{d}{p}}_{p,1}}+||b(\tau)||_{B^{\frac{d}{p}}_{p,1}}\right)d\tau.
\end{eqnarray}
On the other hand, thanks to (\ref{2log-interpolation}), we have
\begin{eqnarray*}
||a(t)||_{B^{\frac{d}{p}}_{p,1}}
&\leq& C ||a(t)||_{B^{\frac{d}{p}}_{p,\infty}}
\ln \left(e+\frac{||a(t)||_{B^{1+\frac{d}{p}}_{p,\infty}}}{||a(t)||_{B^{\frac{d}{p}}_{p,\infty}}}\right)\\ \nonumber
&=& C ||a(t)||_{B^{\frac{d}{p}}_{p,\infty}}
\ln \left(e+\frac{||a(t)||_{B^{1+\frac{d}{p}}_{p,\infty}}\left(1+\frac{||b(t)||_{B^{\frac{d}{p}}_{p,\infty}}}{||a(t)||_{B^{\frac{d}{p}}_{p,\infty}}}\right)}
{||a(t)||_{B^{\frac{d}{p}}_{p,\infty}}+||b(t)||_{B^{\frac{d}{p}}_{p,\infty}}}\right)\\ \nonumber
&\leq& C ||a(t)||_{B^{\frac{d}{p}}_{p,\infty}} \ln \left(e+\frac{K G(t)}{H(t)}\right),
\end{eqnarray*}
and likewise
\begin{eqnarray*}
||b(t)||_{B^{\frac{d}{p}}_{p,1}}
\leq C ||b(t)||_{B^{\frac{d}{p}}_{p,\infty}} \ln \left(e+\frac{K G(t)}{H(t)}\right),
\end{eqnarray*}
which implies that
\begin{eqnarray}\label{4lemma2-proof-a+b-log-estimates}
&&||a(t)||_{B^{\frac{d}{p}}_{p,1}}+||b(t)||_{B^{\frac{d}{p}}_{p,1}}\\ \nonumber
&\leq& C (||a(t)||_{B^{\frac{d}{p}}_{p,\infty}}+||b(t)||_{B^{\frac{d}{p}}_{p,\infty}}) \ln \left(e+\frac{K G(t)}{H(t)}\right).
\end{eqnarray}
Thus, in view of (\ref{4lemma2-assumption}) and the fact that
$$\ln\left(e+\frac{\alpha}{x}\right)\leq \ln(e+\alpha)(1-\ln x),\quad \forall \,\ x\in(0,1],\,\ \alpha>0,$$
we can deduce from (\ref{4lemma2-proof-H-estimates-1}) and (\ref{4lemma2-proof-a+b-log-estimates}) that
\begin{eqnarray*}
H(t)\leq H(0)+C\int_0^t G(\tau)\ln(e+ K G(\tau))H(\tau)(1-\ln H(\tau))d\tau.
\end{eqnarray*}
According to Lemma \ref{2osgood} (set $\mu(r)\triangleq r(1-\ln r)$)
and (\ref{4lemma2-assumption}) again, one has
\begin{eqnarray*}
\frac{H(t)}{e}\leq \left(\frac {H(0)}{e}\right)^{\exp (-C\int_0^t G(\tau) \ln(e+K G(\tau))d\tau)},
\end{eqnarray*}
which leads to the desired result.
In particular, notice that (\ref{4lemma2-initial-assumption}) implies (\ref{4lemma2-assumption}) with $T^\star=T$.
This completes the proof of Lemma 4.2.
\end{proof}

In order to prove Theorem 1.2, it suffices to verify the following lemma.
\begin{lemma4}
Assume that $d\in \mathbb{N_+}$ and  $1\leq p\leq 2d$.
Let $(u_0,\gamma_0)\in B^{1+\frac{d}{p}}_{p,1}(\mathbb{R}^d)\times B^{1+\frac{d}{p}}_{p,1}(\mathbb{R}^d)$
and $(u,\gamma)$ be the corresponding solution to system (\ref{1mchnew}), which is guaranteed by Lemma 4.1.
Then there exist a time $T>0$ and a neighborhood $V$ of $(u_0,\gamma_0)$ in $B^{1+\frac{d}{p}}_{p,1}(\mathbb{R}^d)\times B^{1+\frac{d}{p}}_{p,1}(\mathbb{R}^d)$
such that the mapping $(u_0,\gamma_0)\mapsto (u,\gamma)$ is continuous from $V$ into
$$C([0,T]; B^{1+\frac{d}{p}}_{p,1}(\mathbb{R}^d)\times B^{1+\frac{d}{p}}_{p,1}(\mathbb{R}^d))
\cap C^1{([0,T]; B^{\frac{d}{p}}_{p,1}(\mathbb{R}^d)\times B^{\frac{d}{p}}_{p,1}(\mathbb{R}^d))}.$$
\end{lemma4}

To this end, we need a key convergence result as follows.
\begin{proposition4} \label{4convergence-prop}
\cite{Yan-Yin-DCDS-2015}
Let $d\in \mathbb{N_+}$, $1\leq p\leq \infty$, $1\leq r< \infty$ and $1+\frac{d}{p}<s\neq 2+\frac{d}{p}$
$(\text{or}\,\ s\geq 1+\frac{d}{p}\,\ \text{and}\,\ r=1)$.
Denote $\bar{\mathbb{N}}\triangleq \mathbb{N}\cup \{\infty\}$. Suppose that
$(v^n)_{n\in\bar{\mathbb{N}}}\in C([0,T]; B^{s-1}_{p,r}(\mathbb{R}^d))$ is the solution to
\begin{equation*}
\left\{\begin{array}{ll}
\partial_t v^n +a^n\cdot\nabla v^n=f,\\
{v^n}|_{t=0}=v_0 \end{array}\right.
\end{equation*}
with $v_0\in B^{s-1}_{p,r}(\mathbb{R}^d)$, $f\in L^1(0,T; B^{s-1}_{p,r}(\mathbb{R}^d))$ and that,
for some $\alpha(t)\in L^1(0,T)$ such that
$$\sup\limits_{n\in\bar{\mathbb{N}}} ||\nabla a^n(t)||_{B^{s-1}_{p,r}(\mathbb{R}^d)}\leq \alpha(t).$$
If $a^n$ tends to $a^\infty$ in $L^1(0,T; B^{s-1}_{p,r}(\mathbb{R}^d))$,
then $v^n$ tends to $v^\infty$ in $C([0,T]; B^{s-1}_{p,r}(\mathbb{R}^d))$.
\end{proposition4}

\textbf{Proof of Lemma 4.3.} We divide the proof into three steps as follows.

\textbf{Step 1: Continuity in $C([0,T];B^{\frac{d}{p}}_{p,1}(\mathbb{R}^d)\times B^{\frac{d}{p}}_{p,1}(\mathbb{R}^d)).$}\\
Thanks to Lemmas 4.1-4.2, we have
\begin{eqnarray*}
&& \frac{||v-u||_{L^\infty(0,T; B^{\frac{d}{p}}_{p,\infty})}+||\eta-\gamma||_{L^\infty(0,T; B^{\frac{d}{p}}_{p,\infty})}} {e}\\
&\leq& e^{CMT}\left(\frac{||v_0-u_0||_{B^{\frac{d}{p}}_{p,\infty}}+||\eta_0-\gamma_0||_{B^{\frac{d}{p}}_{p,\infty}}} {e}\right)^{\exp(-CMT\ln(e+K M))},
\end{eqnarray*}
provided that
\begin{eqnarray*}
||v_0-u_0||_{B^{\frac{d}{p}}_{p,\infty}}+||\eta_0-\gamma_0||_{B^{\frac{d}{p}}_{p,\infty}}
\leq e^{1-\exp(CMT\ln(e+K M))}.
\end{eqnarray*}
In view of Lemma 4.1 and (\ref{4interpolation}), we complete the proof of Step 1.

\textbf{Step 2: Continuity in $C([0,T];B^{1+\frac{d}{p}}_{p,1}(\mathbb{R}^d)\times B^{1+\frac{d}{p}}_{p,1}(\mathbb{R}^d)).$}\\
Let $(u^n,\gamma^n)_{n\in\bar{\mathbb{N}}}$ be the solution to the following Cauchy problem:
\begin{equation*}
\left\{\begin{array}{ll}
\partial_t u^n +u^n\cdot\nabla u^n=F_1(u^n,\gamma^n)\triangleq F^n_1(t,x),\\
\partial_t \gamma^n +u^n\cdot\nabla \gamma^n=F_2(u^n,\gamma^n)\triangleq F^n_2(t,x),\\
{u^n}|_{t=0}=u^n_0(x), \\
{\gamma^n}|_{t=0}=\gamma^n_0(x),
\end{array}\right.
\end{equation*}
where $F_1(u^n,\gamma^n), F_2(u^n,\gamma^n)$ are defined by (\ref{1-F1}) and (\ref{1-F2}).

Suppose that $(u^n_0,\gamma^n_0)_{n\in\bar{\mathbb{N}}}\in B^{1+\frac{d}{p}}_{p,1}\times B^{1+\frac{d}{p}}_{p,1}$
and $(u_0^n,\gamma_0^n)$ tends to $(u_0^\infty,\gamma_0^\infty)$ in $B^{1+\frac{d}{p}}_{p,1}\times B^{1+\frac{d}{p}}_{p,1}$.
Thanks to Lemma 4.1, we can find $T,\,\ M>0$ such that for all $n\in \bar{\mathbb{N}}$,
$(u^n,\gamma^n)\in E^{1+\frac{d}{p}}_{p,1}(T)\times E^{1+\frac{d}{p}}_{p,1}(T)$ and
\begin{eqnarray}\label{4lemma3-unibdd}
\sup\limits_{n\in \bar{\mathbb{N}}}
\left(||u^n||_{L^\infty(0,T; B^{1+\frac{d}{p}}_{p,1})}+||\gamma^n||_{L^\infty(0,T; B^{1+\frac{d}{p}}_{p,1})}\right)\leq M.
\end{eqnarray}
According to Step 1,
it suffices to show that
\begin{eqnarray*}
(\nabla u^n,\nabla \gamma^n) \rightarrow (\nabla u^\infty, \nabla \gamma^\infty) \quad \text{in} \quad
C([0,T];B^{\frac{d}{p}}_{p,1}\times B^{\frac{d}{p}}_{p,1}),\,\ \text{as}\,\ n\rightarrow \infty,
\end{eqnarray*}
or in components,
\begin{eqnarray} \label{4lemma3-step2}
\quad (\theta^n_i, \sigma^n) \rightarrow (\theta^\infty_i, \sigma^\infty) \quad \text{in} \quad
C([0,T];B^{\frac{d}{p}}_{p,1}\times B^{\frac{d}{p}}_{p,1}),\,\ \text{as}\,\ n\rightarrow \infty
\end{eqnarray}
with $\theta^n_i\triangleq \nabla(u^n)_{i}$  and  $\sigma^n\triangleq \nabla \gamma^n$ ($n\in \bar{\mathbb{N}}$, $i=1,2,\cdot\cdot\cdot,d$).

Indeed, for each fixed $i\in\{1,2,\cdot\cdot\cdot,d\}$,
$(\theta^n_i, \sigma^n)_{n\in\bar{\mathbb{N}}}$ solves the following transport equations:
\begin{equation*}
\left\{\begin{array}{ll}
\partial_t \theta^n_i +u^n \cdot\nabla \theta^n_i=f^n_i(t,x),\\
\partial_t \sigma^n +u^n \cdot\nabla \sigma^n=g^n(t,x),\\
{\theta^n_i}|_{t=0}=\nabla (u^n_0)_{i},\\
{\sigma^n}|_{t=0}=\nabla \gamma^n_0,
\end{array}\right.
\end{equation*}
where
\begin{eqnarray*}
f^n_i(t,x)\triangleq \nabla(F^n_1)_{i}-\sum\limits_{j=1}^d \partial_j (u^n)_i \nabla {(u^n)_j}
\end{eqnarray*}
and
\begin{eqnarray*}
g^n(t,x)\triangleq \nabla F^n_2-\sum\limits_{j=1}^d \partial_j \gamma^n \nabla {(u^n)_j}.
\end{eqnarray*}

Next, we decompose $\theta^n_i=\theta^{n,1}_i+ \theta^{n,2}_i$
and $\sigma^n=\sigma^{n,1}+\sigma^{n,2}$ for $n\in\mathbb{N}$ with
\begin{equation}\label{4lemma3-step2-eqs1}
\left\{\begin{array}{ll}
\partial_t \theta^{n,1}_i +u^n \cdot\nabla \theta^{n,1}_i=f^n_i-f^\infty_i,\\
\partial_t \sigma^{n,1} +u^n \cdot\nabla \sigma^{n,1}=g^n-g^\infty,\\
{\theta^{n,1}_i}|_{t=0}=\nabla (u^n_0)_{i}-\nabla (u^\infty_0)_{i},\\
{\sigma^{n,1}}|_{t=0}=\nabla \gamma^n_0-\nabla \gamma^\infty_0,
\end{array}\right.
\end{equation}
and
\begin{equation}\label{4lemma3-step2-eqs2}
\left\{\begin{array}{ll}
\partial_t \theta^{n,2}_i +u^n \cdot\nabla \theta^{n,2}_i=f^\infty_i,\\
\partial_t \sigma^{n,2} +u^n \cdot\nabla \sigma^{n,2}=g^\infty,\\
{\theta^{n,2}_i}|_{t=0}=\nabla (u^\infty_0)_{i},\\
{\sigma^{n,2}}|_{t=0}=\nabla \gamma^\infty_0.
\end{array}\right.
\end{equation}
Note that $-(I-\Delta)^{-1}\text{div},\,-(I-\Delta)^{-1}\in Op(S^{-1})$
and $B^{\frac{d}{p}}_{p,1}(\mathbb{R}^d)$ is an algebra.
Thanks to Lemma \ref{2-multiplier} and (\ref{4lemma3-unibdd}),
we can readily gather that $(f^n_i, g^n)_{n\in \bar{\mathbb{N}}}$ is uniformly bounded in
$C([0,T];B^{\frac{d}{p}}_{p,1}\times B^{\frac{d}{p}}_{p,1})$, which implies that, for $i=1,2,\cdot\cdot\cdot,d$,
\begin{eqnarray} \label{4lemma3-step2-RHS}
&& ||f^n_i-f^\infty_i||_{B^{\frac{d}{p}}_{p,1}}+||g^n-g^\infty||_{B^{\frac{d}{p}}_{p,1}}\\ \nonumber
&\leq& CM (||u^n-u^\infty||_{B^{\frac{d}{p}}_{p,1}}+||\nabla u^n-\nabla u^\infty||_{B^{\frac{d}{p}}_{p,1}}
+||\gamma^n-\gamma^\infty||_{B^{\frac{d}{p}}_{p,1}}\\ \nonumber
&&+||\nabla \gamma^n-\nabla \gamma^\infty||_{B^{\frac{d}{p}}_{p,1}}).
\end{eqnarray}
Applying Lemma \ref{2transport-estimates} (i) to system (\ref{4lemma3-step2-eqs1}), one gets
\begin{eqnarray*}
&&||\theta^{n,1}_i(t)||_{B^{\frac{d}{p}}_{p,1}}+||\sigma^{n,1}(t)||_{B^{\frac{d}{p}}_{p,1}}\\
&\leq& e^{C\int_0^t ||\nabla u^n(\tau)||_{B^{\frac{d}{p}}_{p,1}}d\tau}
\{||\nabla (u^n_0)_i-\nabla (u^\infty_0)_i||_{B^{\frac{d}{p}}_{p,1}}+||\nabla \gamma^n_0-\nabla \gamma^\infty_0||_{B^{\frac{d}{p}}_{p,1}}\\
&&+\int_0^t (||(f^n_i-f^\infty_i)(\tau)||_{B^{\frac{d}{p}}_{p,1}}
+||(g^n-g^\infty)(\tau)||_{B^{\frac{d}{p}}_{p,1}})d\tau\},
\end{eqnarray*}
which together with (\ref{4lemma3-unibdd}) and (\ref{4lemma3-step2-RHS}) yield that, for all $t\in [0,T]$,
\begin{eqnarray} \label{4lemma3-step2-endestimates}
&&||\theta^{n,1}_i(t)||_{B^{\frac{d}{p}}_{p,1}}+||\sigma^{n,1}(t)||_{B^{\frac{d}{p}}_{p,1}}\\ \nonumber
&\leq& CM e^{CMT} \{||\nabla u^n_0-\nabla u^\infty_0||_{B^{\frac{d}{p}}_{p,1}}
+||\nabla \gamma^n_0-\nabla \gamma^\infty_0||_{B^{\frac{d}{p}}_{p,1}}\\ \nonumber
&&+\int_0^T (||(u^n-u^\infty)(\tau)||_{B^{\frac{d}{p}}_{p,1}}
+||(\gamma^n-\gamma^\infty)(\tau)||_{B^{\frac{d}{p}}_{p,1}})d\tau\\ \nonumber
&&+\int_0^T (||(\nabla u^n-\nabla u^\infty)(\tau)||_{B^{\frac{d}{p}}_{p,1}}
+||(\nabla\gamma^n-\nabla\gamma^\infty)(\tau)||_{B^{\frac{d}{p}}_{p,1}})d\tau\}.
\end{eqnarray}

On the other hand, notice that $u^n$ tends to $u^\infty$ in $C([0,T];B^{\frac{d}{p}}_{p,1})$ according to Step 1.
Then applying Proposition \ref{4convergence-prop} (with $s=1+\frac{d}{p}$ and $r=1$) to system (\ref{4lemma3-step2-eqs2}),
one infers that $(\theta^{n,2}_i, \sigma^{n,2})$ tends to $(\theta^{\infty,2}_i, \sigma^{\infty,2})$
in $C([0,T];B^{\frac{d}{p}}_{p,1}\times B^{\frac{d}{p}}_{p,1})$.

Thus, for arbitrary $\varepsilon>0$, for $n\in\mathbb{N}$ large enough, by virtue of (\ref{4lemma3-step2-endestimates})
and Step 1, we deduce that for all $t\in [0,T]$,
\begin{eqnarray*}
&&||(\nabla (u^n)_i-\nabla (u^\infty)_i)(t)||_{B^{\frac{d}{p}}_{p,1}}+||(\nabla \gamma^n-\nabla \gamma^\infty)(t)||_{B^{\frac{d}{p}}_{p,1}}\\
&=&||(\theta^n_i-\theta^\infty_i)(t)||_{B^{\frac{d}{p}}_{p,1}}+||(\sigma^n-\sigma^\infty)(t)||_{B^{\frac{d}{p}}_{p,1}}\\
&\leq& ||(\theta^{n,2}_i-\theta^\infty_i)(t)||_{B^{\frac{d}{p}}_{p,1}}+||(\sigma^{n,2}-\sigma^\infty)(t)||_{B^{\frac{d}{p}}_{p,1}}
+||\theta^{n,1}_i(t)||_{B^{\frac{d}{p}}_{p,1}}+||\sigma^{n,1}(t)||_{B^{\frac{d}{p}}_{p,1}}\\
&\leq& \varepsilon+ CM e^{CMT}
\{\varepsilon+||\nabla u^n_0-\nabla u^\infty_0||_{B^{\frac{d}{p}}_{p,1}}+||\nabla \gamma^n_0-\nabla \gamma^\infty_0||_{B^{\frac{d}{p}}_{p,1}}\\
&&+\int_0^T (||(\nabla u^n-\nabla u^\infty)(\tau)||_{B^{\frac{d}{p}}_{p,1}}
+||(\nabla\gamma^n-\nabla\gamma^\infty)(\tau)||_{B^{\frac{d}{p}}_{p,1}})d\tau\},
\end{eqnarray*}
which together with Gronwall's inequality leads to (\ref{4lemma3-step2}).

\textbf{Step 3: Continuity in $C^1([0,T];B^{\frac{d}{p}}_{p,1}(\mathbb{R}^d)\times B^{\frac{d}{p}}_{p,1}(\mathbb{R}^d)).$}\\
According to (\ref{4lemma3-unibdd}) and system (\ref{1mchnew}) itself, Step 1 and Step 2 imply that Step 3 hold true.
Thus, we have proven Lemma 4.3.  Therefore,we complete the proof of Theorem 1.2.
\quad \quad \quad \quad \quad \quad \quad \quad \quad \quad \quad \quad \quad \quad \quad
\quad\quad \quad \quad \quad \quad \quad\quad \quad \quad \quad \quad $\square$

\section{Blow up}
\newtheorem {remark5}{Remark}[section]
\newtheorem{theorem5}{Theorem}[section]
\newtheorem{lemma5}{Lemma}[section]
\newtheorem{proposition5}{Proposition}[section]

In this section, we will prove three blow-up criteria (Theorems 1.3-1.5) of the strong solutions to system (\ref{1mchnew})
by means of the Littlewood-Paley decomposition and the energy method.\\

\textbf{Proof of Theorem 1.3.}
Applying $\Delta_q$ to both sides of the first equation in system (\ref{1mchnew}), one has
\begin{eqnarray}\label{051-eq1}
\partial_t \Delta_q u +u\cdot \nabla (\Delta_q u)
= [u,\Delta_q ] \cdot \nabla u +\Delta_q F_1(u,\gamma).
\end{eqnarray}
Taking the $L^2(\mathbb{R}^d)$ inner product of (\ref{051-eq1}) with $p |\Delta_q u|^{p-2} \Delta_q u$,
integrating by parts and using the H\"{o}lder inequality, one infers that
\begin{eqnarray*}
\frac{d}{d t} ||\Delta_q u||_{L^p(\mathbb{R}^d)}^p
&\leq& p \int_{\mathbb{R}^d} |\Delta_q u|^{p-2} \Delta_q u \left([u,\Delta_q ] \cdot \nabla u +\Delta_q F_1(u,\gamma)\right) d x\\
&&+\int_{\mathbb{R}^d} |\Delta_q u|^p |\nabla u| d x\\
&\leq& ||\nabla u||_{L^\infty(\mathbb{R}^d)} ||\Delta_q u||_{L^p(\mathbb{R}^d)}^p
+ p ||\Delta_q u||_{L^p(\mathbb{R}^d)}^{p-1} ( ||[u,\Delta_q ] \cdot \nabla u||_{L^p(\mathbb{R}^d)} \\
&&+||\Delta_q F_1(u,\gamma)||_{L^p(\mathbb{R}^d)}).
\end{eqnarray*}
Hence,
\begin{eqnarray*}
\frac{d}{d t} ||\Delta_q u||_{L^p(\mathbb{R}^d)}
&\leq& \frac{1}{p} ||\nabla u||_{L^\infty(\mathbb{R}^d)} ||\Delta_q u||_{L^p(\mathbb{R}^d)}
+ ||[u,\Delta_q ] \cdot \nabla u||_{L^p(\mathbb{R}^d)}\\
&&+||\Delta_q F_1(u,\gamma)||_{L^p(\mathbb{R}^d)}.
\end{eqnarray*}
Integrating the above inequality with respect to the time $t$ yields
\begin{eqnarray}\label{051-Delta-u}
\quad\quad ||\Delta_q u||_{L^p(\mathbb{R}^d)}
&\leq& ||\Delta_q u_0||_{L^p(\mathbb{R}^d)}
+\frac{1}{p} \int_0^t ||\nabla u||_{L^\infty(\mathbb{R}^d)} ||\Delta_q u||_{L^p(\mathbb{R}^d)} d\tau\\ \nonumber
&&+\int_0^t  \left(||[u,\Delta_q ] \cdot \nabla u||_{L^p(\mathbb{R}^d)}+||\Delta_q F_1(u,\gamma)||_{L^p(\mathbb{R}^d)}\right) d\tau.
\end{eqnarray}
Multiplying by $2^{q s}$ and taking $l^r$ norm on both sides of (\ref{051-Delta-u}), together with the Minkowski inequality imply
\begin{eqnarray}\label{051-u1}
||u(t)||_{B^s_{p,r}}
&\leq& ||u_0||_{B^s_{p,r}}+\frac{1}{p} \int_0^t ||\nabla u||_{L^\infty} ||u||_{B^s_{p,r}} d\tau
+\int_0^t ||F_1(u,\gamma)||_{B^s_{p,r}} d\tau \nonumber \\
&&+\int_0^t  ||(2^{q s} ||[u,\Delta_q ] \cdot \nabla u||_{L^p})_{q\geq -1}||_{l^r} d\tau.
\end{eqnarray}

Noting that $-(I-\Delta)^{-1}\text{div}, -(I-\Delta)^{-1}\in Op(S^{-1})$
and applying Lemma \ref{2-multiplier} and (\ref{2morse1}), one deduces
\begin{eqnarray}\label{051-F1}
||F_1(u,\gamma)||_{B^s_{p,r}}
&\leq& C(||\nabla u||_{L^\infty} ||\nabla u||_{B^{s-1}_{p,r}}
+||\gamma||_{L^\infty} ||\gamma||_{B^{s-1}_{p,r}}+||\nabla \gamma||_{L^\infty} ||\nabla \gamma||_{B^{s-1}_{p,r}}\nonumber \\
&&+||\nabla u||_{L^\infty} ||u||_{B^{s-1}_{p,r}}+||u||_{L^\infty} ||\nabla u||_{B^{s-1}_{p,r}})\nonumber \\
&\leq& C(||u||_{L^\infty}+||\nabla u||_{L^\infty})||u||_{B^s_{p,r}}+C(||\gamma||_{L^\infty}+||\nabla \gamma||_{L^\infty})||\gamma||_{B^s_{p,r}}.
\end{eqnarray}

Thanks to Lemma \ref{2commutator}, we have
\begin{eqnarray}\label{051-commutator1}
&& ||(2^{q s} ||[u,\Delta_q ] \cdot \nabla u||_{L^p})_{q\geq -1}||_{l^r}\\ \nonumber
&\leq& C \left(||\nabla u||_{L^\infty} ||u||_{B^{s}_{p,r}}+||\nabla u||_{L^\infty} ||\nabla u||_{B^{s-1}_{p,r}}\right)\\ \nonumber
&\leq& C ||\nabla u||_{L^\infty} ||u||_{B^{s}_{p,r}}.
\end{eqnarray}

Thus, by (\ref{051-u1}), (\ref{051-F1}) and (\ref{051-commutator1}), we obtain
\begin{eqnarray}\label{051-u-final}
||u(t)||_{B^s_{p,r}} &\leq& ||u_0||_{B^s_{p,r}}
+C\int_0^t \left(||u||_{L^\infty}+||\nabla u||_{L^\infty}+||\gamma||_{L^\infty}+||\nabla \gamma||_{L^\infty}\right)\nonumber\\
&&\times  \left(||u(\tau)||_{B^s_{p,r}}+||\gamma(\tau)||_{B^s_{p,r}}\right) d\tau.
\end{eqnarray}

On the other hand, applying $\Delta_q$ to the second equation in system (\ref{1mchnew}) yields
\begin{eqnarray}\label{051-eq2}
\partial_t \Delta_q \gamma +u\cdot \nabla (\Delta_q \gamma)
= [u,\Delta_q ] \cdot \nabla \gamma +\Delta_q F_2(u,\gamma).
\end{eqnarray}

Taking the $L^2(\mathbb{R}^d)$ inner product of (\ref{051-eq2}) with $p |\Delta_q \gamma|^{p-2} \Delta_q \gamma$,
and integrating by parts, one gets
\begin{eqnarray*}
\frac{d}{d t} ||\Delta_q \gamma||_{L^p(\mathbb{R}^d)}^p
&\leq& p \int_{\mathbb{R}^d} |\Delta_q \gamma|^{p-2} \Delta_q \gamma \left([u,\Delta_q ] \cdot \nabla \gamma +\Delta_q F_2(u,\gamma)\right) d x\\
&&+\int_{\mathbb{R}^d} |\Delta_q \gamma|^p |\nabla u| d x,
\end{eqnarray*}
which along with the H\"{o}lder inequality leads to
\begin{eqnarray*}
\frac{d}{d t} ||\Delta_q \gamma||_{L^p(\mathbb{R}^d)}
&\leq& \frac{1}{p} ||\nabla u||_{L^\infty(\mathbb{R}^d)} ||\Delta_q \gamma||_{L^p(\mathbb{R}^d)}
+ ||[u,\Delta_q ] \cdot \nabla \gamma||_{L^p(\mathbb{R}^d)}\\
&&+||\Delta_q F_2(u,\gamma)||_{L^p(\mathbb{R}^d)}.
\end{eqnarray*}
Integrating the above inequality with respect to the time $t$ implies
\begin{eqnarray}\label{051-Delta-u}
\quad\quad ||\Delta_q \gamma||_{L^p(\mathbb{R}^d)}
&\leq& ||\Delta_q \gamma_0||_{L^p(\mathbb{R}^d)}
+\frac{1}{p} \int_0^t ||\nabla u||_{L^\infty(\mathbb{R}^d)} ||\Delta_q \gamma||_{L^p(\mathbb{R}^d)} d\tau\\ \nonumber
&&+\int_0^t  \left(||[u,\Delta_q ] \cdot \nabla \gamma||_{L^p(\mathbb{R}^d)}+||\Delta_q F_2(u,\gamma)||_{L^p(\mathbb{R}^d)}\right) d\tau.
\end{eqnarray}
Multiplying by $2^{q s}$ and taking $l^r$ norm on both sides of (\ref{051-Delta-u}), and using the Minkowski inequality, one infers
\begin{eqnarray}\label{051-gamma1}
||\gamma(t)||_{B^s_{p,r}}
&\leq& ||\gamma_0||_{B^s_{p,r}}+\frac{1}{p} \int_0^t ||\nabla u||_{L^\infty} ||\gamma||_{B^s_{p,r}} d\tau
+\int_0^t ||F_2(u,\gamma)||_{B^s_{p,r}} d\tau \nonumber \\
&&+\int_0^t  ||(2^{q s} ||[u,\Delta_q ] \cdot \nabla \gamma||_{L^p})_{q\geq -1}||_{l^r} d\tau.
\end{eqnarray}
Similar to (\ref{051-F1}) and (\ref{051-commutator1}), we can easily get
\begin{eqnarray*}
||F_2(u,\gamma)||_{B^s_{p,r}}
&\leq& C ||\nabla u||_{L^\infty} ||\gamma||_{B^s_{p,r}}+C(||\gamma||_{L^\infty}+||\nabla \gamma||_{L^\infty})||u||_{B^s_{p,r}}
\end{eqnarray*}
and
\begin{eqnarray*}
||(2^{q s} ||[u,\Delta_q ] \cdot \nabla \gamma||_{L^p})_{q\geq -1}||_{l^r}
\leq C \left(||\nabla u||_{L^\infty} ||\gamma||_{B^{s}_{p,r}}+||\nabla \gamma||_{L^\infty} || u||_{B^{s}_{p,r}}\right),
\end{eqnarray*}
which together with (\ref{051-gamma1}) implies
\begin{eqnarray*}\label{051-gamma-final}
||\gamma(t)||_{B^s_{p,r}} &\leq& ||\gamma_0||_{B^s_{p,r}}
+C\int_0^t \left(||\nabla u||_{L^\infty}||\gamma||_{B^s_{p,r}}
+(||\gamma||_{L^\infty}+||\nabla \gamma||_{L^\infty})||u||_{B^s_{p,r}}\right) d\tau.
\end{eqnarray*}
This along with (\ref{051-u-final}) gives
\begin{eqnarray*}
||u(t)||_{B^s_{p,r}} +||\gamma(t)||_{B^s_{p,r}}
&\leq& ||u_0||_{B^s_{p,r}} +||\gamma_0||_{B^s_{p,r}}
+C\int_0^t \left(||u(\tau)||_{B^s_{p,r}}+||\gamma(\tau)||_{B^s_{p,r}}\right)\\
&&\times \left(||u||_{L^\infty}+||\nabla u||_{L^\infty}+||\gamma||_{L^\infty}+||\nabla \gamma||_{L^\infty}\right) d\tau,
\end{eqnarray*}
which together with the Gronwall inequality yields
\begin{eqnarray}\label{051-theorem1-3}
&& ||u(t)||_{B^s_{p,r}} +||\gamma(t)||_{B^s_{p,r}}\\ \nonumber
&\leq& \left(||u_0||_{B^s_{p,r}} +||\gamma_0||_{B^s_{p,r}}\right)
e^{C\int_0^t \left(||u||_{L^\infty}+||\nabla u||_{L^\infty}+||\gamma||_{L^\infty}+||\nabla \gamma||_{L^\infty}\right) d\tau}.
\end{eqnarray}

By virtue of (\ref{051-theorem1-3}) and the Sobolev embedding theorem, we complete the proof of Theorem 1.3.
\quad \quad \quad \quad \quad \quad \quad \quad \quad \quad \quad \quad \quad
\quad\quad \quad \quad \quad \quad \quad\quad \quad \quad \quad \quad $\square$\\

\textbf{Proof of Theorem 1.4.}
Taking the $L^2(\mathbb{R}^d)$ inner product of the first equation in system (\ref{1mchnew}) with $q |u|^{q-2} u$ ($\forall\, q>2$),
integrating by parts and using the H\"{o}lder inequality, one obtains
\begin{eqnarray*}
\frac{d}{d t} ||u||_{L^q(\mathbb{R}^d)}^q
&\leq& \int_{\mathbb{R}^d} |u|^q |\nabla u| d x
+q \int_{\mathbb{R}^d} |u|^{q-2} u F_1(u,\gamma) d x\\
&\leq& ||\nabla u||_{L^\infty(\mathbb{R}^d)} ||u||_{L^q(\mathbb{R}^d)}^q
+ q ||u||_{L^q(\mathbb{R}^d)}^{q-1} ||F_1(u,\gamma)||_{L^q(\mathbb{R}^d)}.
\end{eqnarray*}
Hence,
\begin{eqnarray}\label{052-du}
\frac{d}{d t} ||u||_{L^q(\mathbb{R}^d)}
\leq  ||\nabla u||_{L^\infty(\mathbb{R}^d)} ||u||_{L^q(\mathbb{R}^d)}
+ ||F_1(u,\gamma)||_{L^q(\mathbb{R}^d)}.
\end{eqnarray}

In addition, there exists a constant $c>0$ independent of $q$ such that
\begin{eqnarray}\label{052-F1}
||F_1(u,\gamma)||_{L^q(\mathbb{R}^d)}
&\leq&  c||\nabla u(\nabla u+\nabla u^T)-\nabla u^T \nabla u-\nabla u (\text{div} u)||_{L^q(\mathbb{R}^d)}\nonumber\\
&&+c||\frac{1}{2}(|\nabla u|^2+\gamma^2+|\nabla \gamma|^2)I-\nabla \gamma^T \nabla \gamma||_{L^q(\mathbb{R}^d)}\nonumber\\
&&+c||u(\text{div} u)+u\cdot \nabla u^T||_{L^q(\mathbb{R}^d)}\nonumber\\
&\leq& c ||\nabla u||_{L^\infty(\mathbb{R}^d)} \left(||u||_{L^q(\mathbb{R}^d)}+||\nabla u||_{L^q(\mathbb{R}^d)}\right)\nonumber\\
&&+c \left(||\gamma||_{L^\infty(\mathbb{R}^d)}+||\nabla \gamma||_{L^\infty(\mathbb{R}^d)}\right)
\left(||\gamma||_{L^q(\mathbb{R}^d)}+ ||\nabla \gamma||_{L^q(\mathbb{R}^d)}\right),
\end{eqnarray}
which together with (\ref{052-du}) implies
\begin{eqnarray}\label{052-dudt}
\frac{d}{d t} ||u||_{L^q(\mathbb{R}^d)}
&\leq& c ||\nabla u||_{L^\infty(\mathbb{R}^d)} \left(||u||_{L^q(\mathbb{R}^d)}+||\nabla u||_{L^q(\mathbb{R}^d)}\right)\nonumber\\
&&+c \left(||\gamma||_{L^\infty(\mathbb{R}^d)}+||\nabla \gamma||_{L^\infty(\mathbb{R}^d)}\right)
\left(||\gamma||_{L^q(\mathbb{R}^d)}+ ||\nabla \gamma||_{L^q(\mathbb{R}^d)}\right).
\end{eqnarray}

Applying $\nabla$ to both sides of the first equation in system (\ref{1mchnew}) gives
\begin{eqnarray*}
\partial_t (\nabla u) + \nabla (u\cdot \nabla u) = \nabla F_1(u,\gamma).
\end{eqnarray*}
Taking the $L^2(\mathbb{R}^d)$ inner product of the above equation with $q |\nabla u|^{q-2} \nabla u$ (for any $q>2$),
integrating by parts and using the H\"{o}lder inequality, one obtains
\begin{eqnarray}\label{052-d-du}
\quad \frac{d}{d t} ||\nabla u||_{L^q(\mathbb{R}^d)}
\leq  2||\nabla u||_{L^\infty(\mathbb{R}^d)} ||\nabla u||_{L^q(\mathbb{R}^d)}
+ ||\nabla F_1(u,\gamma)||_{L^q(\mathbb{R}^d)}.
\end{eqnarray}

Similar to (\ref{052-F1}), we have
\begin{eqnarray*}
||\nabla F_1(u,\gamma)||_{L^q(\mathbb{R}^d)}
&\leq& c ||\nabla u||_{L^\infty(\mathbb{R}^d)} \left(||u||_{L^q(\mathbb{R}^d)}+||\nabla u||_{L^q(\mathbb{R}^d)}\right)\nonumber\\
&&+c \left(||\gamma||_{L^\infty(\mathbb{R}^d)}+||\nabla \gamma||_{L^\infty(\mathbb{R}^d)}\right)
\left(||\gamma||_{L^q(\mathbb{R}^d)}+ ||\nabla \gamma||_{L^q(\mathbb{R}^d)}\right),
\end{eqnarray*}
which along with (\ref{052-dudt}) and (\ref{052-d-du}) leads to
\begin{eqnarray*}
&& \frac{d}{d t} \left(||u||_{L^q(\mathbb{R}^d)}+||\nabla u||_{L^q(\mathbb{R}^d)}\right)\\ \nonumber
&\leq& c ||\nabla u||_{L^\infty(\mathbb{R}^d)} \left(||u||_{L^q(\mathbb{R}^d)}+||\nabla u||_{L^q(\mathbb{R}^d)}\right)\\ \nonumber
&&+c \left(||\gamma||_{L^\infty(\mathbb{R}^d)}+||\nabla \gamma||_{L^\infty(\mathbb{R}^d)}\right)
\left(||\gamma||_{L^q(\mathbb{R}^d)}+ ||\nabla \gamma||_{L^q(\mathbb{R}^d)}\right),
\end{eqnarray*}
where $c$ is independent of $q$.\\
Making use of the Gronwall inequality, one gets
\begin{eqnarray*}
||u||_{L^q}+||\nabla u||_{L^q}
&\leq& e^{c \int_0^t ||\nabla u||_{L^\infty} d\tau}
(||u_0||_{L^q}+||\nabla u_0||_{L^q}
+c \int_0^t (||\gamma||_{L^\infty}+||\nabla \gamma||_{L^\infty}) \\
&&\times (||\gamma||_{L^q}+ ||\nabla \gamma||_{L^q}) d\tau ).
\end{eqnarray*}
Letting $q\rightarrow \infty$ and recalling the assumption $B^s_{p,r}(\mathbb{R}^d) \hookrightarrow Lip (\mathbb{R}^d)$, one has
\begin{eqnarray}\label{052-lemma-u-du}
&& ||u||_{L^\infty}+||\nabla u||_{L^\infty}\\ \nonumber
&\leq& C \left(||u_0||_{B^s_{p,r}}+\int_0^t (||\gamma||_{L^\infty}+||\nabla \gamma||_{L^\infty})^2 d\tau \right)
e^{C \int_0^t ||\nabla u||_{L^\infty} d\tau}.
\end{eqnarray}

On the other hand, taking the $L^2(\mathbb{R}^d)$ inner product of the second equation in system (\ref{1mchnew})
with $q |\gamma|^{q-2} \gamma$ ($\forall\, q>2$), integrating by parts and using the H\"{o}lder inequality, one infers
\begin{eqnarray*}
\frac{d}{d t} ||\gamma||_{L^q(\mathbb{R}^d)}^q
&\leq& \int_{\mathbb{R}^d} |\gamma|^q |\nabla u| d x
+q \int_{\mathbb{R}^d} |\gamma|^{q-2} \gamma F_2(u,\gamma) d x\\
&\leq& ||\nabla u||_{L^\infty(\mathbb{R}^d)} ||\gamma||_{L^q(\mathbb{R}^d)}^q
+ q ||\gamma||_{L^q(\mathbb{R}^d)}^{q-1} ||F_2(u,\gamma)||_{L^q(\mathbb{R}^d)}.
\end{eqnarray*}
Or hence,
\begin{eqnarray*}
\frac{d}{d t} ||\gamma||_{L^q(\mathbb{R}^d)}
\leq  ||\nabla u||_{L^\infty(\mathbb{R}^d)} ||\gamma||_{L^q(\mathbb{R}^d)}
+ ||F_2(u,\gamma)||_{L^q(\mathbb{R}^d)}.
\end{eqnarray*}
While similar to the proof of (\ref{052-F1}),
\begin{eqnarray*}
||F_2(u,\gamma)||_{L^q(\mathbb{R}^d)}
\leq  c||\nabla u||_{L^\infty(\mathbb{R}^d)} (||\gamma||_{L^q(\mathbb{R}^d)}+||\nabla\gamma||_{L^q(\mathbb{R}^d)}),
\end{eqnarray*}
where $c$ is independent of $q$.
Then we have
\begin{eqnarray}\label{052-dgammadt}
\frac{d}{d t} ||\gamma||_{L^q(\mathbb{R}^d)}
\leq  c||\nabla u||_{L^\infty(\mathbb{R}^d)} (||\gamma||_{L^q(\mathbb{R}^d)}+||\nabla\gamma||_{L^q(\mathbb{R}^d)}).
\end{eqnarray}

Applying $\nabla$ to both sides of the second equation in system (\ref{1mchnew}) yields
\begin{eqnarray*}
\partial_t (\nabla \gamma) + \nabla (u\cdot \nabla \gamma) = \nabla F_2(u,\gamma).
\end{eqnarray*}
By taking the $L^2(\mathbb{R}^d)$ inner product of the above equation with $q |\nabla \gamma|^{q-2} \nabla \gamma$ ($q>2$),
integrating by parts and using the H\"{o}lder inequality, one gets
\begin{eqnarray}\label{052-d-dgamma}
\quad\quad \frac{d}{d t} ||\nabla \gamma||_{L^q(\mathbb{R}^d)}
&\leq&  \frac{1+q}{q} ||\nabla u||_{L^\infty(\mathbb{R}^d)} ||\nabla \gamma||_{L^q(\mathbb{R}^d)}
+ ||\nabla F_2(u,\gamma)||_{L^q(\mathbb{R}^d)}\\ \nonumber
&\leq&  2||\nabla u||_{L^\infty(\mathbb{R}^d)} ||\nabla \gamma||_{L^q(\mathbb{R}^d)}
+ ||\nabla F_2(u,\gamma)||_{L^q(\mathbb{R}^d)}.
\end{eqnarray}

Similar to (\ref{052-F1}), we have
\begin{eqnarray*}
||\nabla F_2(u,\gamma)||_{L^q(\mathbb{R}^d)}
\leq c||\nabla u||_{L^\infty(\mathbb{R}^d)} (||\gamma||_{L^q(\mathbb{R}^d)}+||\nabla\gamma||_{L^q(\mathbb{R}^d)}),
\end{eqnarray*}
which along with (\ref{052-dgammadt}) and (\ref{052-d-dgamma}) ensures
\begin{eqnarray*}
\frac{d}{d t} \left(||\gamma||_{L^q(\mathbb{R}^d)}+||\nabla \gamma||_{L^q(\mathbb{R}^d)}\right)
\leq c ||\nabla u||_{L^\infty(\mathbb{R}^d)} \left(||\gamma||_{L^q(\mathbb{R}^d)}+||\nabla \gamma||_{L^q(\mathbb{R}^d)}\right),
\end{eqnarray*}
where $c$ is independent of $q$.\\
Thanks to the Gronwall inequality again, we have
\begin{eqnarray*}
||\gamma||_{L^q}+||\nabla \gamma||_{L^q}
\leq \left(||\gamma_0||_{L^q}+||\nabla \gamma_0||_{L^q}\right) e^{c \int_0^t ||\nabla u||_{L^\infty} d\tau}.
\end{eqnarray*}
Letting $q\rightarrow \infty$, one infers
\begin{eqnarray}\label{052-lemma-gamma-dgamma}
||\gamma||_{L^\infty}+||\nabla \gamma||_{L^\infty}
&\leq& \left(||\gamma_0||_{L^\infty}+||\nabla \gamma_0||_{L^\infty}\right) e^{c \int_0^t ||\nabla u||_{L^\infty} d\tau}\\ \nonumber
&\leq& C ||\gamma_0||_{B^s_{p,r}} e^{C \int_0^t ||\nabla u||_{L^\infty} d\tau}.
\end{eqnarray}

Combining (\ref{052-lemma-u-du}) with (\ref{052-lemma-gamma-dgamma}), we obtain
\begin{eqnarray*}
&& ||u||_{L^\infty}+||\nabla u||_{L^\infty}+||\gamma||_{L^\infty}+||\nabla \gamma||_{L^\infty}\\
&\leq& C\left(||u_0||_{B^s_{p,r}}+||\gamma_0||_{B^s_{p,r}}+||\gamma_0||^2_{B^s_{p,r}} t\right)
e^{C \int_0^t ||\nabla u||_{L^\infty} d\tau},
\end{eqnarray*}
which together with (\ref{051-theorem1-3}) and the Sobolev embedding theorem complete the proof of Theorem 1.4.
\quad \quad \quad \quad \quad \quad \quad \quad \quad \quad \quad \quad
\quad\quad \quad \quad \quad \quad \quad\quad \quad \quad \quad \quad $\square$ \\

\textbf{Proof of Theorem 1.5.}
In view of Remark {\ref{1-maximal-time}}, we may assume $s>3+\frac{d}{p}$ to prove the theorem.
By taking the $L^2(\mathbb{R}^d)$ inner product of the first equation in system (\ref{1mch}) with $q |m|^{q-2} m$ ($\forall\, q>d$),
integrating by parts and using the H\"{o}lder inequality, we have
\begin{eqnarray*}
\frac{d}{d t} ||m||_{L^q}^q
&\leq& (1+2q)\int_{\mathbb{R}^d} |m|^q |\nabla u| d x
+q \int_{\mathbb{R}^d} |m|^{q-2} m \rho \nabla \gamma d x\\
&\leq& (1+2q) ||\nabla u||_{L^\infty} ||m||_{L^q}^q
+ q ||m||_{L^q}^{q-1} ||\rho||_{L^q} ||\nabla \gamma||_{L^\infty}.
\end{eqnarray*}
Hence,
\begin{eqnarray}\label{053-dm}
\quad\quad \frac{d}{d t} ||m||_{L^q}
&\leq&  \frac{1+2q}{q} ||\nabla u||_{L^\infty} ||m||_{L^q}
+ ||\nabla \gamma||_{L^\infty} ||\rho||_{L^q}\\ \nonumber
&\leq&  3||\nabla u||_{L^\infty} ||m||_{L^q}
+ ||\nabla \gamma||_{L^\infty} ||\rho||_{L^q}.
\end{eqnarray}

On the other hand, taking the $L^2(\mathbb{R}^d)$ inner product of the second equation in system (\ref{1mch})
with $q |\rho|^{q-2} \rho$ ($\forall\, q>d$), integrating by parts and using the H\"{o}lder inequality, one obtains
\begin{eqnarray*}
\frac{d}{d t} ||\rho||_{L^q}
\leq \frac{1+q}{q} ||\nabla u||_{L^\infty} ||\rho||_{L^q}
\leq  2 ||\nabla u||_{L^\infty} ||\rho||_{L^q},
\end{eqnarray*}
which together with (\ref{053-dm}) yields
\begin{eqnarray}\label{053-dm+drho}
\quad\quad \frac{d}{d t} (||m||_{L^q}+||\rho||_{L^q})
\leq  3 (||\nabla u||_{L^\infty}+||\nabla \gamma||_{L^\infty})
(||m||_{L^q}+ ||\rho||_{L^q}).
\end{eqnarray}

Note that $u=(I-\Delta)^{-1} m$ and $\gamma=(I-\Delta)^{-1} \rho$ imply
\begin{eqnarray}\label{053-La-bdd1}
||D^k u||_{L^a(\mathbb{R}^d)}\leq c ||m||_{L^a(\mathbb{R}^d)}
\end{eqnarray}
and
\begin{eqnarray}\label{053-La-bdd2}
||D^k \gamma||_{L^a(\mathbb{R}^d)}\leq c ||\rho||_{L^a(\mathbb{R}^d)},
\end{eqnarray}
where $1\leq a\leq \infty$, $k=0,1,2$, and the constant $c$ is independent of $a$.

In view of (\ref{053-La-bdd1}) and (\ref{053-La-bdd2}), thanks to (\ref{2new-log-interpolation}), we infer that
\begin{eqnarray} \label{053-log-u}
||\nabla u||_{L^\infty(\mathbb{R}^d)}
&\leq& c \frac{2q-d}{q-d}
\left(1+||\nabla u||_{B^{0}_{\infty,\infty}(\mathbb{R}^d)} \ln(e+||\nabla u||_{W^{1,q}(\mathbb{R}^d)})\right) \nonumber\\
&\leq& c \frac{2q-d}{q-d}
\left(1+||\nabla u||_{B^{0}_{\infty,\infty}(\mathbb{R}^d)} \ln(e+||m||_{L^q (\mathbb{R}^d)})\right).
\end{eqnarray}
Likewise,
\begin{eqnarray} \label{053-log-gamma}
\quad\quad ||\nabla \gamma||_{L^\infty(\mathbb{R}^d)}
\leq c \frac{2q-d}{q-d}
\left(1+||\nabla \gamma||_{B^{0}_{\infty,\infty}(\mathbb{R}^d)} \ln(e+||\rho||_{L^q (\mathbb{R}^d)})\right).
\end{eqnarray}

Substituting (\ref{053-log-u}) and (\ref{053-log-gamma}) into (\ref{053-dm+drho}), one gets
\begin{eqnarray*}
\frac{d}{d t} (||m||_{L^q}+||\rho||_{L^q})
&\leq& c \frac{2q-d}{q-d} \left(1+||\nabla u||_{B^{0}_{\infty,\infty}}+||\nabla \gamma||_{B^{0}_{\infty,\infty}} \right)\\
&&\times (||m||_{L^q}+||\rho||_{L^q}) \ln \left(e+(||m||_{L^q}+||\rho||_{L^q})\right).
\end{eqnarray*}
Integrating the above inequality with respect to the time $t$ gives
\begin{eqnarray}\label{053-m+rho-1}
&& ||m(t)||_{L^q}+||\rho(t)||_{L^q}\\ \nonumber
&\leq& ||m_0||_{L^q}+||\rho_0||_{L^q}
+ c \frac{2q-d}{q-d} \int_0^t \left(1+||\nabla u||_{B^{0}_{\infty,\infty}}+||\nabla \gamma||_{B^{0}_{\infty,\infty}} \right)\\ \nonumber
&&\times (||m(\tau)||_{L^q}+||\rho(\tau)||_{L^q}) \ln \left(e+(||m(\tau)||_{L^q}+||\rho(\tau)||_{L^q})\right) d\tau,
\end{eqnarray}
where $c$ is independent of $q$.

Let $q\rightarrow \infty$ in (\ref{053-m+rho-1}). Then we have
\begin{eqnarray}\label{053-m+rho-final}
&& ||m(t)||_{L^\infty}+||\rho(t)||_{L^\infty}\\ \nonumber
&\leq& ||m_0||_{L^\infty}+||\rho_0||_{L^\infty}
+ c \int_0^t \left(1+||\nabla u||_{B^{0}_{\infty,\infty}}+||\nabla \gamma||_{B^{0}_{\infty,\infty}} \right)\\ \nonumber
&&\times (||m(\tau)||_{L^\infty}+||\rho(\tau)||_{L^\infty}) \ln \left(e+(||m(\tau)||_{L^\infty}+||\rho(\tau)||_{L^\infty})\right) d\tau.
\end{eqnarray}

Set $\Phi(t)\triangleq e+||m(t,\cdot)||_{L^\infty}+||\rho(t,\cdot)||_{L^\infty}$.
From (\ref{053-m+rho-final}), one has
\begin{eqnarray}\label{053-phi-1}
\quad\quad\quad\quad  \Phi(t) \leq \Phi(0)
+C\int_0^t \left(1+||\nabla u||_{B^{0}_{\infty,\infty}}+||\nabla \gamma||_{B^{0}_{\infty,\infty}} \right) \Phi(\tau) \ln\Phi(\tau) d\tau.
\end{eqnarray}

Applying Lemma \ref{2osgood} (set $\mu(r)\triangleq r\ln r (r\geq e)$) to (\ref{053-phi-1}) yields
\begin{eqnarray*}
\ln (\ln \Phi(t))
\leq \ln (\ln \Phi(0))
+C\int_0^t \left(1+||\nabla u||_{B^{0}_{\infty,\infty}}+||\nabla \gamma||_{B^{0}_{\infty,\infty}} \right)d\tau,
\end{eqnarray*}
or hence,
\begin{eqnarray} \label{053-phi-final}
\Phi(t) \leq (\Phi(0))^{\exp\left(C\int_0^t (1+||\nabla u||_{B^{0}_{\infty,\infty}}+||\nabla \gamma||_{B^{0}_{\infty,\infty}}) d\tau\right)}.
\end{eqnarray}

According to (\ref{053-La-bdd1}), (\ref{053-La-bdd2}) and (\ref{053-phi-final}), we deduce
\begin{eqnarray*}
&& ||u(t)||_{L^\infty}+||\nabla u(t)||_{L^\infty}+||\gamma(t)||_{L^\infty}+||\nabla \gamma(t)||_{L^\infty}\\
&\leq& C(||m(t)||_{L^\infty}+||\rho(t)||_{L^\infty})\\
&\leq& C (e+||m_0||_{L^\infty}+||\rho_0||_{L^\infty})^{\exp\left(C\int_0^t (1+||\nabla u||_{B^{0}_{\infty,\infty}}+||\nabla \gamma||_{B^{0}_{\infty,\infty}}) d\tau\right)},
\end{eqnarray*}
which along with Theorem 1.3 completes the proof of Theorem 1.5.
\quad\quad \quad \quad $\square$\\

\section{Appendix}

In this Appendix, we give the details that how to rewrite system (\ref{1mch}) to its nonlocal form system (\ref{1mchnew}).
For this, we first introduce some notations:\\
Let $u=(u_1,u_2,\cdot\cdot\cdot,u_d)$, $v=(v_1,v_2,\cdot\cdot\cdot,v_d)$ be vector fields, and
$A=(a_{ij})_{d\times d}$, $B=(b_{ij})_{d\times d}$ be $d\times d$ matrices. Then
$$(i)\quad u\cdot \nabla v\triangleq \sum\limits_{j=1}^d u_j\partial_j v
=u(\nabla v)^T=u\nabla v^T,$$
where $A^T$ denotes the transpose of $A$.
$$(ii)\quad \text{div}u\triangleq \sum\limits_{j=1}^d \partial_j u_j,\,\ \,\ \text{while}\,\ \,\
\text{div}A \triangleq (\text{div}A_1,\text{div}A_2,\cdot\cdot\cdot,\text{div}A_d)$$
with $A=\left(\begin{array}{c}
             A_1\\
           \vdots\\
           A_d\\
          \end{array}\right)$
and each component $A_j=(a_{j1},a_{j2},\cdot\cdot\cdot,a_{jd})$.
In particular, we have $\text{div}(\nabla u)=\Delta u=(\Delta u_1,\Delta u_2,\cdot\cdot\cdot, \Delta u_d)$.
$$(iii)\quad  A:B\triangleq \sum\limits_{i,j=1}^d a_{ij} b_{ij}\,\ \,\ \text{and}\,\ \,\
|A|\triangleq (A:A)^{1/2}.$$

Set $\gamma\triangleq \bar{\rho}-\bar{\rho}_0$.
From the first equation in system (\ref{1mch}), we deduce
\begin{eqnarray}\label{6-1}
&&(I-\Delta) (\partial_t u+u\cdot \nabla u)\\ \nonumber
&=& \partial_t m+u\cdot \nabla u-\Delta (u\cdot \nabla u)\\ \nonumber
&=& u\cdot \nabla (\Delta u)-\Delta(u\cdot \nabla u)+\nabla u^T\cdot(\Delta u)+(\Delta u)\text{div}u
-u\cdot \nabla u^T-u(\text{div}u)\\ \nonumber
&& -\gamma \nabla \gamma+ (\Delta\gamma)\nabla\gamma.
\end{eqnarray}
While
\begin{eqnarray}\label{6-2}
\quad\quad  u\cdot \nabla (\Delta u)-\Delta(u\cdot \nabla u)
=-\text{div}(\nabla u\nabla u+\nabla u\nabla u^T)+(\nabla u)\cdot \nabla(\text{div}u),
\end{eqnarray}

\begin{eqnarray}\label{6-3}
\nabla u^T\cdot(\Delta u)
&=&\text{div}(\nabla u^T \nabla u)-\frac{1}{2}\nabla (|\nabla u|^2)\\ \nonumber
&=&\text{div}\left(\nabla u^T \nabla u-\frac{1}{2}|\nabla u|^2 I \right),
\end{eqnarray}

\begin{eqnarray}\label{6-4}
(\nabla u)\cdot \nabla(\text{div}u)+(\Delta u)\text{div}u
=\text{div}(\nabla u(\text{div}u)),
\end{eqnarray}
and
\begin{eqnarray}\label{6-5}
\gamma \nabla \gamma- (\Delta\gamma)\nabla\gamma
&=&\nabla \left(\frac{1}{2}\gamma^2+\frac{1}{2}|\nabla \gamma|^2\right)-\text{div}(\nabla \gamma^T \nabla \gamma)\\ \nonumber
&=&\text{div}\left(\frac{1}{2}(\gamma^2+|\nabla \gamma|^2) I -\nabla \gamma^T \nabla \gamma\right).
\end{eqnarray}
So, in view of (\ref{6-1})-(\ref{6-5}), one gets the first equation in system (\ref{1mchnew}).

On the other hand, from the second equation in system (\ref{1mch}), we have
\begin{eqnarray}\label{6-6}
&&(I-\Delta) (\partial_t \gamma+u\cdot \nabla \gamma)\\ \nonumber
&=& \partial_t \rho+u\cdot \nabla\gamma-\Delta (u\cdot \nabla \gamma)\\ \nonumber
&=& u\cdot \nabla (\Delta \gamma)-\Delta(u\cdot \nabla \gamma)+(\Delta \gamma)(\text{div}u)-\gamma(\text{div}u).
\end{eqnarray}
While
\begin{eqnarray}\label{6-7}
\quad\quad  u\cdot \nabla(\Delta \gamma)-\Delta(u\cdot \nabla\gamma)
=-\text{div}\left( \nabla\gamma \nabla u +(\nabla\gamma)\cdot\nabla u \right)+ \nabla\gamma \cdot \nabla(\text{div}u),
\end{eqnarray}
and
\begin{eqnarray}\label{6-8}
(\Delta \gamma)(\text{div}u)
=\text{div}(\nabla\gamma (\text{div}u))-\nabla \gamma\cdot \nabla(\text{div}u).
\end{eqnarray}
Hence, by (\ref{6-6})-(\ref{6-8}), we obtain the second equation in system (\ref{1mchnew}).
\quad\quad \quad \quad $\square$

\bigskip
\noindent\textbf{Acknowledgments}
\ The author thanks the referees for their valuable comments and suggestions.
He was partially supported by NNSFC (No.11501226).

\end{document}